\documentclass[11pt,a4paper]{article}
\usepackage{amsmath,amssymb,amsthm}
\usepackage{color}
\usepackage{ascmac}
\usepackage{dsfont}

\usepackage{enumerate}

\theoremstyle{definition}
 \newtheorem{dfn}{Definition}[section]
 
 \newtheorem{lemma}[dfn]{Lemma}
  \newtheorem{rmk}[dfn]{Remark}

\theoremstyle{plain}
 \newtheorem{thm}[dfn]{Theorem}
 \newtheorem{prop}[dfn]{Proposition}
 \newtheorem{lem}[dfn]{Lemma}
 \newtheorem{cor}[dfn]{Corollary}

\newcommand{\bu}{{\mathbf u}}
\newcommand{\bz}{{\mathbf z}}
\newcommand{\bv}{{\mathbf v}}
\newcommand{\bw}{{\mathbf w}}

\newcommand{\bB}{{\mathbf B}}

\newcommand{\bF}{{\mathbf F}}
\newcommand{\bG}{{\mathbf G}}

\newcommand{\bJ}{{\mathbf J}}

\newcommand{\bU}{{\mathbf U}}

\newcommand{\dv}{{\rm div}\,}

\newcommand{\BR}{{\mathbb R}}

\newcommand{\BN}{{\mathbb N}}

\newcommand{\Z}{{\mathbb Z}}

\newcommand{\CB}{{\mathcal B}}

\newcommand{\CD}{{\mathcal D}}

\newcommand{\CF}{{\mathcal F}}
\newcommand{\CI}{{\mathcal I}}

\newcommand{\CL}{{\mathcal L}}

\newcommand{\CN}{{\mathcal N}}

\newcommand{\CS}{{\mathcal S}}

\newcommand{\CP}{{\mathcal P}}

\newcommand{\CV}{{\mathcal V}}

\newcommand{\fh}{{\mathfrak h}}

\newcommand{\fp}{{\mathfrak p}}

\newcommand{\be}{{\mathbf e}}
\newcommand{\bg}{{\mathbf g}}
\newcommand{\bh}{{\mathbf h}}
\newcommand{\pd}{\partial}
\newcommand{\curl}{{\rm curl\,}}

\usepackage{xcolor}

\newcommand{\ep}{\varepsilon}

\newcommand{\wt}{\widetilde}
\newcommand{\wh}{\widehat}

\newcommand{\R}{{\mathbb R}}

\newcommand{\Bf}{{\mathbf f}}
\newcommand{\Bg}{{\mathbf g}}

\newcommand{\hB}{\dot{B}}
\newcommand{\hF}{\dot{F}} 

\newcommand{\hR}{\dot{R}}
\newcommand{\hS}{\dot{S}}
\newcommand{\hT}{\dot{T}}

\newcommand{\hDelta}{\dot{\Delta}}

\newcommand{\ddj}{\dot\Delta_j}
\newcommand{\ddu}{\dot\delta\bu}
\newcommand{\ddB}{\dot\delta\bB}
\newcommand{\ddv}{\dot\delta\bv}

\numberwithin{equation}{section} 

\usepackage{scalerel,stackengine}
\stackMath
\newcommand\reallywidehat[1]{%
\savestack{\tmpbox}{\stretchto{%
  \scaleto{%
    \scalerel*[\widthof{\ensuremath{#1}}]{\kern-.6pt\bigwedge\kern-.6pt}%
    {\rule[-\textheight/2]{1ex}{\textheight}}
  }{\textheight}%
}{0.5ex}}%
\stackon[1pt]{#1}{\tmpbox}%
}

\newcommand{\vertiii}[1]{{\left\vert\kern-0.25ex\left\vert\kern-0.25ex\left\vert #1 
    \right\vert\kern-0.25ex\right\vert\kern-0.25ex\right\vert}}

\begin{document}

\title{\bf On steady solutions of the Hall-MHD system in Besov spaces}
\author{Jin Tan
\thanks{Laboratoire de Math\'{e}matiques AGM, UMR CNRS 8088, Cergy Paris Universit\'{e}, 2 Avenue Adolphe Chauvin, 95302, Cergy-Pontoise Cedex, France\endgraf 
e-mail address:jin.tan@cyu.fr
}
and Hiroyuki Tsurumi
\thanks{Graduate School of Technology, Industrial and Social Sciences, Tokushima University, Tokushima 770--8506, Japan. \endgraf 
e-mail address: tsurumi.hiroyuki@tokushima-u.ac.jp}
and Xin Zhang \thanks{School of Mathematical Sciences,
Tongji University, 
No.1239, Siping Road, Shanghai (200092), China. \endgraf
e-mail address: xinzhang2020@tongji.edu.cn}}

\maketitle

\begin{abstract}
In this paper,  we investigate the well-posedness and ill-posedness issues for 
the incompressible stationary Hall-magnetohydrodynamic (Hall-MHD) system in $\BR^3.$
We first show the existence and uniqueness of solutions provided with the forces in $\hB^{3/p-3}_{p,r}(\R^3)$ for $1\leq p <3$ and $r=1$. 
Moreover, this result can be extended to any $1\leq r\leq \infty$ whenever $p=2,$ without any additional assumption on the physical parameters. 
On the other hand, we establish some ill-posedness results for Hall-MHD system by using the
discontinuity of the solution mapping of the three-dimensional stationary Navier-Stokes equations in \emph{critical}  function spaces $\hB^{3/p-1}_{p,r}(\R^3)$ ($p\geq 3$).\\

\vskip1pc\noindent
keywords: Hall-magnetohydrodynamic system, well-posedness, ill-posedness, Besov spaces
exterior domain
\end{abstract}
%

\section{Introduction}

\subsection{Model and results}
In this paper, we consider the stationary viscous 
Hall-magnetohydrodynamics (Hall-MHD) system  in $\BR^3$ as follows
\begin{equation}\label{eq:hMHD_1}
	\left\{ \begin{aligned}
&\bu\cdot\nabla \bu  + \nabla \fp -\mu \Delta \bu
-(\nabla \times \bB)\times \bB=\Bf, \\
&\dv \bu=0,\\
&-\nu \Delta \bB -\nabla\times  \big((\bu-h \nabla \times \bB)
\times \bB \big) =\bg. \\
\end{aligned}
\right.
\end{equation}
Here $\bu=(u_1(x), u_2(x), u_3(x))^\perp$, $\bB=(B_1(x), B_2(x), B_3(x))^\perp$ and $\fp=\fp(x)$ with $x\in\R^3$ denote the velocity field, the magnetic field, and the pressure, respectively, while $\Bf=(f_1(x), f_2(x), f_3(x))^\perp$ and $\bg=(g_1(x), g_2(x), g_3(x))^\perp$ denote the given external forces.
In addition, $\mu>0$, $\nu>0$, and $h>0$ are given parameters denoting the fluid viscosity, the magnetic resistivity, and the Hall-coefficient, respectively.
For the detailed explanations of this system as a physical model, we refer to \cite{ADFL2011,Hu2003,SU1997} for example.

Actually, there are many previous studies on the non-stationary problem of the Hall-MHD system.
For example, Chae-Degond-Liu \cite{CDL2014} showed the local well-posedness for the general initial data in $H^s$ ($s<5/2$).
After that, Chae-Lee \cite{CL2014} and Ye \cite{Ye2015} obtained the blow-up criterion and the global well-posedness for small data.
Then Wan-Zhou \cite{WZ2015} and Benvenutti-Ferreira \cite{BF2016} considered the initial data allowing low regularity.
In scaling critical functional spaces, Danchin-Tan \cite{DT2021,DT2022} showed the well-posedness in Sobolev and Besov spaces, and Nakasato \cite{Na2022} in Fourier-Besov spaces.

On the other hand, compared to the non-stationary system, there seems to be few studies on the stationary system \eqref{eq:hMHD_1}.
Recently, Li-Su \cite{LS2022} showed the Liouville type theorems for \eqref{eq:hMHD_1} with $\Bf=\bg=0$ and Cho-Neustupa-Yang \cite{CNY2024} generalized it. 
However, the problem on the existence and uniqueness of solutions of \eqref{eq:hMHD_1}  in $\BR^3$ with given external forces has not been fully studied.
In fact, such problem has been well considered in the stationary Navier-Stokes equations (we can regard this as \eqref{eq:hMHD_1} in the case $\bg\equiv\bB\equiv0$).
For the well-posedness results of the Navier-Stokes problem, we refer to the works of Chen \cite{Ch1993} (Sobolev space), Kozono-Yamazaki \cite{KY1995} (Morrey space), and Kaneko-Kozono-Shimizu \cite{KKS} (Besov space). 
On the contrary, Tsurumi \cite{Ts2019a, Ts2019} and Li-Yu-Zhu \cite{LYZ2022} showed the ill-posedness in critical Besov spaces.
Therefore, our motivation is to study the well-posedness and ill-posedness problem of \eqref{eq:hMHD_1}, and to compare the result to that of Navier-Stokes equations.
\medskip

For this purpose, this paper studies \eqref{eq:hMHD_1} in some critical functional spaces.
Let us first explain what we mean by \emph{critical} regularity for the stationary Hall-MHD system \eqref{eq:hMHD_1}. In fact, because the Hall term
$$ h \nabla\times  \big( (\nabla \times \bB)
\times \bB \big) $$
 breaks the natural scaling of the classical  stationary MHD system (i.e. taking $h=0$ in \eqref{eq:hMHD_1}), there does not exists a genuine scaling invariance structure. However,
similar to  the non-stationary Hall-MHD problem in \cite{DT2021},  if we introduce an auxiliary  vector field $\bJ:=\nabla\times \bB$ (that can be interpreted as electron current in physics), then the  following extended system in terms of   $(\bu, \bB, \bJ)$  
\begin{equation}\label{eq:HMHD-ext}
\left\{\begin{aligned}
&-\mu\Delta \bu+\nabla \fp
-(\nabla\times \bB)\times \bB=\Bf,\\
&\dv \bu=0,\\
& -\nu\Delta \bB-\nabla\times \big( (\bu -h\bJ) \times \bB)=\bg,\\
 & -\nu\Delta \bJ-\nabla\times\bigl(\nabla\times\big( (\bu -h \bJ) 
 \times {\rm{curl}}^{-1} \bJ \big)\bigr)=\nabla\times \bg
\end{aligned}\right.
\end{equation}
does have the scaling invariant structure:
\begin{equation}\label{Hall-scaling}
(\bu , \bB , \bJ )(x)\leadsto  \lambda (\bu, \bB, \bJ)( \lambda x),\quad
 \fp(x) \leadsto  \lambda^2 \fp ( \lambda x),\quad 
 (\Bf , \bg , \nabla \times \bg )(x)\leadsto 
 \lambda^3 (\Bf , \bg , \nabla \times \bg )(\lambda x).
\end{equation} 
Above, we have used the so-called \emph{Biot-Savart}  operator
$$\curl^{-1}:=(-\Delta)^{-1} \nabla \times.$$
In view of \eqref{Hall-scaling}, the velocity field $\bu$ in \eqref{eq:HMHD-ext} has the same scaling invariance as the stationary incompressible Navier-Stokes equations.  
So we have the following results on the well-posedness issue of \eqref{eq:hMHD_1}.
\begin{thm}\label{thm:DS_1}
Let $1\leq p<3$ and $1\leq r\leq 2.$
Then \eqref{eq:hMHD_1} is well-posed
\footnote{For the rigorous definition of the well-posedness or the ill-posedness, see the subsection \ref{subsec:WI}. } 
from $D_1$ to $S_1$ with
\begin{equation*}
\begin{aligned}
D_1&:=\big\{(\Bf,\bg): \Bf,\nabla \times \bg \in  \hB^{3/p-3}_{p,1}(\R^3)^3, \,
  \bg\in \hB^{3/p-3}_{p,r}(\R^3)^3, \dv \bg=0 \big\},\\
S_1&:=\big\{(\bu,\bB): \bu, \nabla \times \bB\in \hB^{3/p-1}_{p,1} (\R^3)^3, \,
 \bB \in  \hB^{3/p-1}_{p,r} (\R^3)^3,  \dv\bu=\dv \bB=0 \big\}.
\end{aligned}
\end{equation*}
\end{thm}

Compared with the well-posedness result on the stationary Navier-Stokes equations in \cite{KKS}, the well-posedness issue of \eqref{eq:hMHD_1} for the indices 
$(p,r)\in [1,3)\times (2,\infty]$ is not clear because of the additional Hall term in \eqref{eq:hMHD_1}.
However, using a cancellation property, namely
\footnote{The symbol $(\cdot,\cdot)_{L^2}$ denotes the inner product in $L^2(\R^3).$},
\begin{equation}\label{Cancellation}
\Big(\nabla\times \big( (\nabla\times \bB)\times \bB \big), \bB\Big)_{L^2}
=\big( (\nabla\times \bB)\times \bB, \nabla\times \bB\big)_{L^2}=0,
\end{equation}
 we can prove the existence and the uniqueness of the solution of \eqref{eq:hMHD_1} when $p=2$ and $r\in [1,\infty]$, even though the continuity of the solution mapping might fail. 
\begin{thm} \label{thm:DS_2}
For any $1\leq r\leq \infty,$ we define
\begin{equation*}
\begin{aligned}
D_2&:=\big\{(\Bf,\bg): \Bf,\bg,\nabla \times 
\bg \in  \hB^{-3/2}_{2,r}(\R^3)^3,\dv \bg=0 \big\},\\
S_2&:=\big\{(\bu,\bB): \bu,  \bB, \nabla \times \bB\in \hB^{1/2}_{2,r} (\R^3)^3,
\dv \bu=\dv\bB=0\big\}.
\end{aligned}
\end{equation*}
There exists a positive constant $\delta$ depending only
 of $\mu, \nu,h$ such that if
\begin{equation}\label{small-con-p=2}
 \|(\Bf, \bg, \nabla \times \bg)\|_{\hB^{-3/2}_{2,r}(\R^3)}< \delta 
\end{equation}
for any $(\Bf, \bg)\in D_2,$ then system \eqref{eq:hMHD_1} admits a unique solution $(\bu, \bB)\in S_2.$ 
\end{thm}

Finally, we consider the ill-posedness issue of \eqref{eq:hMHD_1} in the case $p\geq 3,$
which can be derived from the ill-posedness of the stationary Navier-Stokes equations
\begin{equation}\label{eq:NS_1}
	\left\{ \begin{aligned}
&\bu\cdot\nabla \bu  + \nabla \fp -\mu \Delta \bu
=\Bf, \\
&\dv \bu=0. \\
\end{aligned}
\right.
\end{equation}
For \eqref{eq:NS_1}, we can prove the following statement together with some partial results in \cite{Ts2019a, Ts2019, LYZ2022}.
\begin{thm} \label{thm:DS_3}
For any $3\leq p\leq\infty$ and $1\leq r\leq \infty$, we define
\footnote{See Section \ref{sec:besov} for the definition of the homogeneous Besov spaces. }
\begin{equation*}
D_3:=\hB^{3/p-3}_{p,r}(\R^3)^3, \quad
S_3:=\big\{\bu\in \hB^{3/p-1}_{p,r} (\R^3)^3,\dv \bu=0\big\}.
\end{equation*}
\begin{enumerate}[$(1)$]	
\item If $p=3$ and $3/2< r\leq 2$, then \eqref{eq:NS_1} is well-posed from $D_3$ to $S_3$.
\item  If $3<p\leq\infty$ and $1\leq r\leq \infty$, or if $p=3$ and $r\in [1,3/2)\cup (2,\infty]$, then \eqref{eq:NS_1} is ill-posed from $D_3$ to $S_3$
by the discontinuity of solution mapping near the origin.
\end{enumerate}
\end{thm}
Let us list some comments on Theorem \ref{thm:DS_3}.
\begin{itemize}
\item Actually, the result (1) in Theorem \ref{thm:DS_3} for $r=2$ has already shown in \cite{LYZ2022}. In this paper, we extend the well-posedness result of \eqref{eq:NS_1} to the situation where $3/2<r<2.$ 

\item  The result (2) of Theorem \ref{thm:DS_3} in the case 
\begin{equation*}
(p,r)\in  (3,\infty]\times [1,\infty]  \,\,\,\text{or}\,\,\,
p=3\,\,\,\text{and}\,\,\,
2< r\leq \infty
\end{equation*}
has shown by the second author \cite{Ts2019a, Ts2019}. The main aim of Theorem \ref{thm:DS_3} (2) extends the ill-posedness result to the case $p=3$ and $1\leq r<3/2.$ In other words,
for the continuity of the solution mapping of \eqref{eq:NS_1} from $\hB^{-2}_{3,r}(\R^3)^3$ to $\hB^{0}_{3,r}(\R^3)^3,$ both indices $r=3/2$ and $r=2$ are the borderline cases between the ill-posedness issue and the well-posedness issue.
In contrast to the four-dimensional stationary Navier-Stokes problem, $r=2$ is the only borderline case for the resulting solution mapping from $\hB^{-2}_{4,r}(\R^4)^4$ to $\hB^{-2}_{4,r}(\R^4)^4$ which is observed in \cite{LYZ2022}. 
\end{itemize}
\medskip

Actually, we can show Theorem \ref{thm:DS_3} (2) by constructing a sequence $\{\Bf_n\}_{n\in\mathbb{N}}\subset D_3$ such that
\begin{enumerate}[(i)]
\item For each $\Bf=\Bf_n$, there exists a solution $\bu=\bu_n\in S_3$ of  \eqref{eq:NS_1};
\item $\|\Bf_n\|_{D_3}\to 0$ as $n\to\infty$, however, $\|\bu_n\|_{S_3}\nrightarrow 0$ as $n\to\infty$. Moreover, this sequence $\{\bu_n\}_{n\in\mathbb{N}}$ has an additional property as follows:
\begin{enumerate}[(a)]
\item If $p=\infty$ and $1\leq r \leq \infty$, then there exists a constant $c_0>0$ such that $\|\bu_n\|_{\hB^{-1}_{\infty,\infty}}>c_0$ for every $n\in\mathbb{N}$ (see \cite{Ts2019a}). We note here that $\hB^{-1}_{\infty,\infty}$ is the weakest critical space including $S_3$.
\item If $(p,r)\in  (3,\infty)\times [1,\infty]  \,\,\,\text{or}\,\,\,
p=3\,\,\,\text{and}\,\,\,
2< r\leq \infty$, then we also see that $\|\bu_n\|_{\hB^{-1}_{\infty,\infty}}\nrightarrow 0$ as $n\to\infty$ (see \cite{Ts2019}). 
\item If $p=3$ and $1\leq r <3/2$, then there exists a constant $c_0>0$ such that $\|\bu_n\|_{S_3}>c_0$ for every $n\in\mathbb{N}.$ This is one of new results of this paper (see Section \ref{sec:ill} for more details).

\end{enumerate}
\end{enumerate}
As we may regard $(\bu_n, 0)$ as a sequence of the solution of \eqref{eq:hMHD_1} for 
$(\Bf, \bg)=(\Bf_n,0)$,  then it is immediate to see the following result.
\begin{cor}
Let $3<p\leq\infty$ and $1\leq r\leq \infty$, or let $p=3$ and $r\in [1,3/2)\cup (2,\infty].$
The problem \eqref{eq:hMHD_1} is ill-posed from $D_4$ to $S_4$ with
\begin{equation*}
\begin{aligned}
D_4&:=\big\{(\Bf,\bg): \Bf \in  \hB^{3/p-3}_{p,r}(\R^3)^3, \,
  \bg\in \hB^{0}_{1,1}(\R^3)^3, \dv \bg=0 \big\},\\
S_4&:=\big\{(\bu,\bB): \bu\in \hB^{3/p-1}_{p,r} (\R^3)^3, \,
 \bB \in  \hB^{-1}_{\infty,\infty} (\R^3)^3,  \dv\bu=\dv \bB=0 \big\}.
\end{aligned}
\end{equation*}

\end{cor}

\subsection{Main idea}
In this subsection, we shall explain the idea to prove Theorems \ref{thm:DS_1}, \ref{thm:DS_2}
and \ref{thm:DS_3}.

\begin{itemize}

\item The key to the proof of  Theorem \ref{thm:DS_1} is to consider an appropriate reformulation of the system \eqref{eq:HMHD-ext}, suitably rewritten in the form of a generalized  stationary Navier-Stokes system.
Recall the algebraic identities as presented in \cite{DT2021}:
\begin{equation*}
\bz\cdot\nabla \bw = \dv (\bw\otimes \bz),\quad 
 (\nabla\times \bw)\times  \bw
 =( \bw\cdot\nabla)\bw-\nabla(|\bw|^{2}/2),
\end{equation*}
where $\bigl({\rm{div}}(\bw\otimes \bz)\bigr)^{j}{:=}\sum_{k=1}^{3}\partial_{k}(w^{j}z^{k}).$
Notice the fact that $\dv \bB=0$ for $\dv \bg=0.$
Then applying the Leray projector $\CP:=\CI+\nabla(-\Delta)^{-1}\dv$ to \eqref{eq:HMHD-ext} yields that
\begin{equation}\label{eq:hMHD_3}
	\left\{ \begin{aligned}
&-\mu \Delta \bu = \CP \Bf-\CP \dv  (\bu \otimes \bu) 
+ \CP \dv (\bB\otimes \bB), \\
&-\nu \Delta \bB =\CP \bg
+\nabla\times\big((\bu-h \bJ)\times \bB \big),  \\
&-\nu \Delta \bJ 
=\CP\bh+\nabla \times \Big( \nabla\times  
\big((\bu-h \bJ)\times \curl^{-1}\bJ  \big)  \Big), 
\end{aligned}
\right.
\end{equation}
with $\bh:=\nabla \times \bg.$ 
It is not hard to construct the solution of \eqref{eq:hMHD_3} by the fixed point argument.

\item The proof of Theorem \ref{thm:DS_2} relies on the cancellation property \eqref{Cancellation}, which was suggested in Danchin-Tan \cite{DT2021} in the case that $\mu=\nu.$ More precisely, in the case that $\mu=\nu$, by considering the auxiliary vector-field  $\bv:=\bu-h\bJ$ (that can be interpreted as velocity of electron), the only quasilinear term cancels out when performing an energy method, since
\begin{equation}\label{Cancellation2}
\Big( \nabla\times\big( (\nabla\times \bv)\times \bB\big), \bv\Big)_{L^2}=0.
\end{equation}

In the stationary setting for \emph{general} $\mu, \nu,$ the auxiliary vector-field $\bv$ satisfies
\begin{equation}\label{eq-v}
\left\{
\begin{aligned}
&  -\mu\Delta \bu=\bB\cdot\nabla \bB-\bu\cdot\nabla \bu-\nabla \fp+\Bf,\\
&  -\nu\Delta \bB=\nabla\times(\bv\times \bB)+\bg,\\
& -\nu\Delta \bv=\bB\cdot\nabla \bB-\bu\cdot\nabla \bu
-h\nabla\times \big( (\nabla\times \bv)\times \bB \big)+\Bf-h\nabla\times \bg\\
&\hspace{2cm}+\nabla\times (\bv\times \bu)
+2\nabla\times \big( \bv\cdot\nabla \curl^{-1}(\bu-\bv)\big)
-\nabla \fp+(\mu-\nu)\Delta \bu,\\
&\dv \bu=\dv \bB=\dv \bv=0.
\end{aligned}
\right.
\end{equation}
After localization of the above system by means of the Littlewood-Paley  spectral 
cut-off operators $\dot \Delta_j$ (defined in the next section), the  
identity  \eqref{Cancellation2} still holds, up to some lower order manageable commutator term (see Proposition \ref{Prop-commu}). 
In order to  control the term $(\mu-\nu) \Delta \bu$ in $\eqref{eq-v}_3$ when $\mu\neq \nu,$  we perform a renormalized energy method.
This enables us to prove the existence and uniqueness results for $(\Bf, \bg, \Bf-h\nabla\times \bg)$ in the critical space $\dot B^{-3/2}_{2,r}$ ($1\leq r \leq \infty$).

\item As for Theorem \ref{thm:DS_3}, we first check the ill-posedness result when $p>3$ by direct application of the method in \cite{Ts2019}.
Secondly, we show the well-posedness when $p=3$ and $3/2< r\leq 2$ by using the embedding relation between homogeneous Besov spaces and Triebel-Lizorkin spaces, and applying the H\"{o}lder inequality. 
After that, using this new well-posedness result, we prove that the sequence of external forces proposed by \cite{LYZ2022} may cause the ill-posedness even when $p=3$ with $1\leq r<3/2$.

\end{itemize}

\section{Toolbox in analysis}
\label{sec:besov}
In this section, we introduce the definitions and some important properties concerning  homogeneous Besov spaces. Then we shall give some abstract theory 
on the well-posedness and ill-posedness due to \cite{BT2006}
\subsection{Functional spaces}
First, we introduce the Littlewood-Paley decomposition (see \cite{BCD2011,Sa2018} for instance).
Consider some smooth radial function $\varphi$ supported in the annulus $\{ \xi \in \R^3 : {3}/{4} \leq |\xi| \leq {8}/{3} \}$ satisfying the following decomposition
\begin{align} \label{eq:L-P_decomp}
 \sum_{j \in \mathbb{Z} } \varphi(2^{-j}\xi) =1,
 ~~~ \forall ~\xi \in \R^3 \backslash \{0\}.
\end{align}
Then we define  Fourier truncation operators
\begin{equation*}
\hDelta_j := \varphi(2^{-j}D),\quad 
\hS_j := \sum_{j'\leq j-1}\hDelta_{j'},\,\,\forall j\in \Z.
\end{equation*}

\begin{dfn} 
Let $\CP(\R^3)$ and $\CS'(\R^3)$ denote the polynomial space and the tempered distribution space in $\R^3.$
For $s\in\R$ and $(p,r) \in [1,\infty]^2,$ the homogeneous Besov space  $\hB^s_{p,r}(\R^3)$ is defined by
\begin{equation}\label{def:besov}
\hB^{s}_{p,r}(\R^3):= \big\{ \,u \in \mathcal{S}'(\R^3)/\CP(\R^3): \|u\|_{\hB^s_{p,r}} :=\big\|2^{js}\|\hDelta_j u \|_{L^p(\R^3)} \big\|_{\ell^r(\mathbb{Z}) } <\infty \,\big\}.
\end{equation}
Sometimes, we also write $\hB^s_{p,r}$ for short.
\end{dfn}

One advantage of homogeneous Besov space is the following property
\begin{equation}\label{eq:homo}
\|f(\lambda \cdot )\|_{\hB^{s}_{p,r}(\R^3)} \sim 
\lambda^{s-3\slash p}\|f(\cdot )\|_{\hB^{s}_{p,r}(\R^3)} 
\,\,\, \text{for any}\,\,\lambda>0.
\end{equation}
Here $A \sim B$ stands for $C_1 A \leq B \leq C_2 A$ for some constants $C_1$ and $C_2.$ 
Moreover, the following result is well-known,
\begin{equation}\label{eq:embd_Besov}
\hB^{s}_{p_1,r_1}(\R^3) \hookrightarrow 
\hB^{s-3(1\slash p_1-1\slash p_2)}_{p_2,r_2}(\R^3).
\end{equation}
for any $s\in \R,$ $1\leq p_1\leq p_2 \leq \infty$ and $1\leq r_1 \leq r_2\leq \infty.$  
\medskip

To handle the nonlinear terms in \eqref{eq:hMHD_1} in Besov spaces setting, 
we introduce the following \emph{paraproduct} and \emph{remainder} operators (after J.-M. Bony in \cite{Bo1981}):
$$\hT_u v := \sum_{j\in \Z} \hS_{j-1}u \hDelta_j v\quad\hbox{and}\quad 
\hR(u,v) := \sum_{\genfrac{}{}{0pt}{}{j\in\Z}{|j-k|\leq 1}} \hDelta_{j} u \hDelta_{k} v.
$$
Then  any product may be formally decomposed as follows:
\begin{equation}\label{eq:bony}
uv = \hT_u v + \hT_v u + \hR(u,v).
\end{equation} 
Thanks to \eqref{eq:bony} and the continuity of paraproduct and the remainder operators (see \cite[Section 2.6]{BCD2011} for instance), it is not hard to prove the following proposition.
For simplicity, we use $A \lesssim B$ for $A \leq C B$ up to some harmless constant $C$ hereafter.
\begin{prop}\label{prop:BE_1}
Let $1\leq r,r_1,r_2 \leq \infty.$ Then the following assertions hold true. 
\begin{enumerate}
\item If $1\leq p<6,$ then we have
\begin{equation}\label{es:pl_2}
\|uv\|_{\hB^{3/p-1}_{p,r}(\R^3)} 
\lesssim \|u\|_{\hB^{3/p-1}_{p,r}(\R^3)}
\|v\|_{\hB^{3/p}_{p,\infty}(\R^3) \cap L^{\infty}(\R^3)}.
\end{equation}

\item If $1\leq p<3$ and $1/r_1 +1/r_2=1/r,$ then we have
\begin{equation}\label{es:pl_1}
\|uv\|_{\hB^{3/p-2}_{p,r}(\R^3)} 
\lesssim \|u\|_{\hB^{3/p-1}_{p,r_1}(\R^3)}
\|v\|_{\hB^{3/p-1}_{p,r_2}(\R^3)}.
\end{equation}
\end{enumerate}
\end{prop}
\begin{proof}
We first prove \eqref{es:pl_2}.
By the decomposition \eqref{eq:bony}  and the continuity of paraproduct operator (see \cite[Section 2.6]{BCD2011} for instance),
we have 
\begin{align*}
\|\hT_u v \|_{\hB^{3/p-1}_{p,r}} 
&\lesssim \|u\|_{\hB^{-1}_{\infty,r}}
\|v\|_{\hB^{3/p}_{p,\infty}} 
\lesssim \|u\|_{\hB^{3/p-1}_{p,r}}
\|v\|_{\hB^{3/p}_{p,\infty}},\\
 \|\hT_v u \|_{\hB^{3/p-1}_{p,r}} 
&  \lesssim \|v\|_{L^{\infty}}
\|u\|_{\hB^{3/p-1}_{p,r}}.
\end{align*}
For the remainder part, if $1\leq p<3,$ that is, $3/p-1>0,$ then we have 
\begin{equation*}
\|\hR(u,v) \|_{\hB^{3/p-1}_{p,r}} 
\lesssim \|u\|_{\hB^{3/p-1}_{p,r}}
\|v\|_{\hB_{\infty,\infty}^0} 
\lesssim \|u\|_{\hB^{3/p-1}_{p,r}}
\|v\|_{\hB^{3/p}_{p,\infty}} .
\end{equation*}
If $3\leq p<6,$ then we obtain that 
\begin{equation*}
 \|\hR(u,v) \|_{\hB^{3/p-1}_{p,r}} \lesssim 
\|\hR(u,v) \|_{\hB^{6/p-1}_{p/2,r}} 
\lesssim \|u\|_{\hB^{3/p-1}_{p,r}}
\|v\|_{\hB^{3/p}_{p,\infty}}. 
\end{equation*}

Next, notice that
\begin{align*}
\|\hT_u v \|_{\hB^{3/p-2}_{p,r}} 
&  \lesssim \|u\|_{\hB^{-1}_{\infty,r_1}}
\|v\|_{\hB^{3/p-1}_{p,r_2}} 
\lesssim \|u\|_{\hB^{3/p-1}_{p,r_1}}
\|v\|_{\hB^{3/p-1}_{p,r_2}},\\
 \|\hT_v u \|_{\hB^{3/p-2}_{p,r}} 
&  \lesssim \|v\|_{\hB^{-1}_{\infty,r_2}}
\|u\|_{\hB^{3/p-1}_{p,r_1}} 
\lesssim\|v\|_{\hB^{3/p-1}_{p,r_2}} 
 \|u\|_{\hB^{3/p-1}_{p,r_1}}.
\end{align*}
For the remainder term,  we first suppose $1\leq p <3/2,$  i.e. $3/p-2>0,$ and then we have
\begin{equation*}
 \|\hR(u,v) \|_{\hB^{3/p-2}_{p,r}} 
\lesssim \|u\|_{\hB^{-1}_{\infty,r_1}}
\|v\|_{\hB^{3/p-1}_{p,r_2}}
\lesssim \|u\|_{\hB^{3/p-1}_{p,r_1}}
\|v\|_{\hB^{3/p-1}_{p,r_2}}. 
\end{equation*}
Otherwise, if $3/2\leq p<3,$ then we suppose that $3/p=1+\sigma$ for some $0<\sigma\leq 1.$ Let $q_{\sigma}:=6/(4+\sigma).$ Then we have 
\begin{equation*}
1<q_\sigma <3/2, \quad  
\frac{1}{q_{\sigma}} = \frac{1}{p} + \frac{1}{q}
\end{equation*}
with $q=6/(2-\sigma) \in (3,6].$ Thus we have 
\begin{equation*}
 \|\hR(u,v) \|_{\hB^{3/p-2}_{p,r}} 
\lesssim  \|\hR(u,v) \|_{\hB^{q_\sigma-2}_{q_{\sigma},r}} 
\lesssim \|u\|_{\hB^{3/p-1}_{p,r_1}}
\|v\|_{\hB^{3/q-1}_{q,r_2}} 
\lesssim \|u\|_{\hB^{3/p-1}_{p,r_1}}
\|v\|_{\hB^{3/p-1}_{p,r_2}}.
\end{equation*}
Thus we have \eqref{es:pl_1}.
\end{proof}
\medskip

Next, we establish some commutator estimate in order to study \eqref{eq:hMHD_1} in $L^2$ type spaces for Theorem \ref{thm:DS_2}.
\begin{prop} \label{Prop-commu}
Let $s\in (-3/2,3/2].$ 
For all $r, \rho_2\in[1, \infty]$  and    $\rho_1 \in(2, \infty],$ we have 
\begin{equation}\label{com2}
\bigl\|2^{js}\| [\dot\Delta_j, b]a\|_{L^2(\R^3)}\bigr\|_{\ell^r(\Z)}
\lesssim \|b\|_{ \dot B^{2/\rho_1+3/2}_{2, \infty}(\R^3)} 
\|a\|_{\dot B^{s-2/\rho_1}_{2,r}(\R^3)}
+\|b\|_{\dot B^{s+3/2+2/\rho_2}_{2,r}(\R^3)}
\|a\|_{\dot B^{-2/\rho_2}_{2, \infty}(\R^3)},
\end{equation}
where we have set
\begin{align*}
[\ddj, b]a:=\ddj (b\,a)-b\,\ddj a.
\end{align*}
\end{prop}

\begin{proof}
Proving \eqref{com2} relies  on the decomposition 
\begin{equation}\label{decompo}
[\dot\Delta_j, b]a=[\dot\Delta_j,\hT_b]a
+ \dot\Delta_j \big( \hT_ab+\hR(a,b)\big)
-\big( \hT_{\dot\Delta_ja}b+\hR(\dot\Delta_ja,b) \big).
\end{equation}
By the localization properties of the Littlewood-Paley decomposition,
the first term of \eqref{decompo} may be decomposed into
$$[\ddj,\hT_b]a =\sum_{|j'-j|\leq4} [\ddj, \dot S_{j'-1}b]\dot\Delta_{j'}a.$$
Now, according to \cite[Lemma 2.97]{BCD2011},  we have 
$$\| [\ddj, \dot S_{j'-1}b]\dot\Delta_{j'}a\|_{L^2}\lesssim 2^{-j}\|\nabla\dot S_{j'-1}b\|_{L^\infty}\|\dot\Delta_{j'}a\|_{L^2}.$$
Since $2/\rho_1-1<0$, 
$$\|\nabla S_{j'-1}b\|_{ L^\infty}\lesssim 2^{j'(1-2/\rho_1)} 
\|\nabla b\|_{ \dot B^{2/\rho_1-1}_{\infty,\infty}}.$$
Hence, for all $(j,j')\in\Z^2,$  we have
$$2^{js} \|[\ddj, \dot S_{j'-1}b]\dot\Delta_{j'}a\|_{ L^2}\lesssim 
2^{(j-j')(s-1)}2^{j'(s-2/\rho_1)}\|\dot\Delta_{j'}a\|_{ L^2}
\|\nabla b\|_{ \dot B^{2/\rho_1-1}_{\infty,\infty}}.$$
Therefore,  summing up on $j'$ fulfilling $|j-j'|\leq 4,$ then taking the $\ell^r(\Z)$ norm,
$$\bigl\|2^{js}\|[\ddj,\hT_b]a  \|_{ L^2}\bigr\|_{\ell^r} 
\lesssim  \|\nabla b\|_{ \dot B^{2/\rho_1-1}_{\infty,\infty}} 
\|a\|_{ \dot B^{s-2/\rho_1}_{2,r}}.$$
\medskip

The next two terms in \eqref{decompo} may be bounded according
to  \cite[Theorems 2.47, 2.52]{BCD2011}
$$\bigl\|2^{js}\|\ddj \hT_ab \|_{ L^2}\|_{\ell^r} 
+\bigl\|2^{js}\|\ddj \hR(a,b)\|_{ L^2}\|_{\ell^r} 
\lesssim \|a\|_{ \dot B^{-2/\rho_2}_{2, \infty}}
\|  b\|_{\ \dot B^{s+3/2+2/\rho_2}_{2,r}}$$
 for any $s+3/2>0$ and $1\leq \rho_2\leq  \infty.$
\medskip

Finally, owing to the properties of localization of the Littlewood-Paley decomposition,  one has
\begin{equation*} 
\hT_{\dot\Delta_ja}b+R(\dot\Delta_ja,b)=\sum_{j'\geq j-2} \dot S_{j'+2}\dot\Delta_ja\,\dot\Delta_{j'}b.
\end{equation*}
Then we take advantage of the fact that 
$$\| \dot S_{j'+2}\ddj a\|_{ L^\infty}\lesssim \|\ddj a\|_{ L^\infty} 
\lesssim  2^{(3/2+2/\rho_2)j}
\|a\|_{ \dot B^{-3/2-2/\rho_2}_{\infty,\infty}},$$ 
which implies that
 $$\begin{aligned}
 2^{js}\| \hT_{\ddj a}b+R(\ddj a,b)\|_{ L^2}
&\lesssim \sum_{j'\geq j-2}2^{js}\|\dot S_{j'+2}a\|_{ L^\infty}\|\dot\Delta_{j'}b\|_{ L^2}\\
 &\lesssim  \|a\|_{ \dot{B}^{-3/2-2/\rho_2}_{\infty, \infty}}\!\sum_{j'\geq j-2} \!\!
 2^{(s+3/2+2/\rho_2)(j-j')}\,2^{(s+3/2+2/\rho_2)j'}\!\|\dot\Delta_{j'} b\|_{L^2}.
\end{aligned}$$
 Taking the $\ell^r(\Z)$ norm of both sides and using Young's inequality,  
 we end up with 
 $$
 \bigl\|2^{js}\| \hT_{\ddj a}b+R(\ddj a,b)\|_{L^2}\bigr\|_{\ell^r(\Z)} \lesssim
  \|a\|_{ \dot B^{-2/\rho_2}_{2,\infty}}
\| b\|_{ \dot B^{s+3/2+2/\rho_2}_{2,r}},$$
where we have used that  $s+3/2+2/\rho_2>0.$
This completes the proof of \eqref{com2}. 
\end{proof}

\subsection{Abstract theory on the well-posedness issue}
\label{subsec:WI}
According to \cite{BT2006}, we review some abstract theory on the well-posedness and ill-posedness of the following abstract equation
\begin{equation}\tag{${A}$}\label{eq:A}
u= Lf + N(u,u).
\end{equation}
Here $f\in D$ is a given data, and $u\in S$ is an unknown solution of \eqref{eq:A}
for some Banach spaces $(D,\|\cdot\|_D)$ and $(S,\|\cdot\|_S).$ 
In addition, $L: D \to S$ is a densely defined linear operator, and $N: S\times S\to S$ denotes a densely defined bilinear form.  

\begin{dfn}\label{def:wp}
We call that the equation \eqref{eq:A} is well-posed from $(D,\|\cdot\|_D)$ to $(S,\|\cdot\|_S)$ if there exist two constants $\varepsilon, \delta>0$ such that
\footnote{Here  $B_X(r)$ denotes the open ball in Banach space $X$ centred at the origin with the radius $r>0.$}
\begin{enumerate}[(i)]
\item for every $f\in B_D(\varepsilon)$, there exist a solution $u\in B_S(\delta)$ of \eqref{eq:A};

\item  if there exist two solutions $u_1, u_2 \in B_S(\delta)$ of \eqref{eq:A} for the same datum $f\in B_D(\varepsilon),$ then $u_1= u_2$ in $S$;

\item The solution map $f \in (B_D(\varepsilon), \|\cdot\|_D) \mapsto u\in (B_S(\delta), \|\cdot\|_S)$, well-defined by {\rm (i)} and {\rm (ii)}, is continuous.
\end{enumerate}
In addition, \eqref{eq:A} is ill-posed from $D$ to $S$ if \eqref{eq:A} is not well-posed from $D$ to $S$.
\end{dfn}
In Definition \ref{def:wp} above, (i) and (ii) state the existence and uniqueness of solutions for sufficiently small data respectively,  while (iii) is about the continuously dependence of solutions on given data. 
On the other hand, we introduce the so-called \emph{quantitatively} well-posedness for \eqref{eq:A} as follows. 
\begin{dfn}\label{def:qwp}
The abstract problem \eqref{eq:A} is called quantitatively well-posed from $(D,\|\cdot\|_D)$ to $(S,\|\cdot\|_S)$ if there exist positive constants $C_1$ and $C_2$ depending solely on $D$ and $S$ such that
\begin{equation}\label{cdt:L}
\|Lf\|_S\leq C_1 \|f\|_D\hspace{5pt} { for\ any\ } f\in D,
\end{equation}
\begin{equation}\label{cdt:N_0}
\|N(u,v)\|_S \leq C_2 \|u\|_S \|v\|_S
\hspace{5pt} { for\ any\ } u, v\in S.
\end{equation}
\end{dfn}

In fact, the quantitatively well-posedness gives us not only the well-posedness in Definition \ref{def:wp}, but also the continuity for some approximation scheme due to \cite[Theorem 3 and Proposition 1]{BT2006} as follows.
\begin{prop}\label{prop:BT}
Suppose that \eqref{eq:A} is quantitatively well-posed from $(D, \|\cdot\|_D)$ to $(S, \|\cdot\|_S)$. 

\begin{enumerate}[$(1)$]
\item \eqref{eq:A} is well-posed from $(D, \|\cdot\|_D)$ to $(S, \|\cdot\|_S)$ in the sense of Definition \ref{def:wp}.  More precisely, if we define the nonlinear maps $A_m: D\to S$ for $m\in\mathbb{N}$ by
\begin{equation}
\begin{cases}
A_1 f := Lf,\\
A_m f := \mathop{\sum}\limits_{k,l\geq 1, k+l=m} N(A_{k} f, A_{l} f),\ \ m\geq 2,
\end{cases}
\nonumber
\end{equation}
then there exists some positive constant $C$ independent of $m$ such that
\begin{equation}\label{es:A_m}
\|A_m f\|_{S} \leq \left(C \|f\|_{D}\right)^m,\,\, \forall\,\, m\in\mathbb{N}.
\end{equation}

Moreover, set that
\begin{equation}
\label{u_f}
u(f):=\sum_{m=1}^\infty A_m f \hspace{10pt}{ in}\ S. 
\end{equation}
Then there exist some constants $C_0$ and (sufficiently small) $\ep$ such that \eqref{eq:A} admits a unique solution $u=u(f)$ in $B_S(C_0\ep)$ for every $f \in B_D(\varepsilon).$

\item Suppose that $D$ and $S$ are given other weaker norms $\|\cdot\|_{\tilde D}$ and $\|\cdot\|_{\tilde S}$ respectively. That is,
\[ \|f\|_{\tilde D}\leq C\|f\|_{D},\hspace{10pt} \|u\|_{\tilde S}\leq C\|u\|_{S} \]
for some constant $C>0$. Assume that the solution map $f \mapsto u$ of \eqref{eq:A} is continuous from $(B_D(\varepsilon), \|\cdot\|_{\tilde D})$ to $(B_S(\delta), \|\cdot\|_{\tilde S})$. Then for every $m\in\mathbb{N}$, $A_m :D\to S$ is also continuous from $(B_D(\varepsilon), \|\cdot\|_{\tilde D})$ to $(B_S(\delta), \|\cdot\|_{\tilde S})$.
\end{enumerate}
\end{prop}

\section{Well-posedness issue in critical spaces}
In this section, we prove Theorem \ref{thm:DS_1} by constructing the solution of \eqref{eq:hMHD_3}. For simplicity, we introduce 
\begin{equation*}
\bU:=\begin{pmatrix}
\bu\\ \bB\\\bJ
\end{pmatrix},\,\,\,
\bF:=\begin{pmatrix}
\Bf\\ \bg\\ \bh
\end{pmatrix}
\,\,\, \text{and}\,\,\,\,
\CL\bF:=\begin{pmatrix}
(-\mu\Delta)^{-1}\CP\Bf\\
 (-\nu\Delta)^{-1}\CP\bg\\  
 (-\nu\Delta)^{-1}\CP\bh
\end{pmatrix}.
\end{equation*}
Then the system \eqref{eq:hMHD_3} can be rewritten as
\begin{equation}\label{eq:hMHD}
\bU= \CL\bF + \CN(\bU,\bU)
\end{equation}
with 
\begin{equation}\label{def:N}
\begin{aligned}
\CN(\bU,\bU)&:=\big(\CN_1(\bU,\bU),\CN_2(\bU,\bU),\CN_3(\bU,\bU)\big),\\
\CN_1(\bU,\bU)&:=(-\mu\Delta)^{-1} \CP \dv( -\bu \otimes \bu+\bB\otimes \bB),\\
\CN_2(\bU,\bU)&:= (-\nu\Delta)^{-1}\nabla\times  \big((\bu-h \bJ)\times \bB \big),\\
\CN_3(\bU,\bU)&:= (-\nu \Delta)^{-1}  \nabla \times \Big( \nabla\times  
\big((\bu-h \bJ)\times \curl^{-1}\bJ  \big)  \Big).
\end{aligned}
\end{equation}

Now, we shall establish the following result for \eqref{eq:hMHD}, where $\bJ$ and $\bh$ are not necessarily dependent of $\bB$ and $\bg$ respectively.
\begin{thm}\label{thm:U_1} 
Let $1\leq p<3$ and $1\leq r\leq 2.$ 
Then the system \eqref{eq:hMHD} is well-posed from $\CD_1$ to $\CS_1$ where
\footnote{For simplicity, we also use the symbol $\hB^{s}_{p,r}(\R^3)$ for vector-valued functions.} 
\begin{equation*}
\begin{aligned}
\CD_1 &:=\big\{(\Bf,\bg,\bh): \Bf,\bh \in \hB^{3/p-3}_{p,1}(\R^3), \,
 \bg \in  \hB^{3/p-3}_{p,r}(\R^3)\big\},\\
\CS_1&:=\big\{(\bu,\bB,\bJ): \bu, \bJ \in \hB^{3/p-1}_{p,1}(\R^3), \,
\bB\in \hB^{3/p-1}_{p,r}(\R^3),\dv \bu=\dv \bB=\dv \bJ=0\big\}.
\end{aligned}
\end{equation*}
In particular, \eqref{eq:hMHD} is well-posed from $\hB^{3/p-3}_{p,1}(\R^3)$ to $\hB^{3/p-1}_{p,1}(\R^3).$ 
\end{thm}
It is enough to prove Theorem \ref{thm:DS_1} by applying Theorem \ref{thm:U_1} with $\bh=\nabla \times \bg.$ So we only prove the well-posedness of \eqref{eq:hMHD} in the following.
According to Proposition \ref{prop:BT}, we shall verify that \eqref{eq:hMHD} is quantitatively well-posed from $\CD_1$ to $\CS_1.$ Indeed, we shall see that
\begin{align}
\|\CL\bF\|_{\CS_1}&\leq C \|\bF\|_{\CD_1}, \label{eq:L1}\\
\|\CN(\bU,\bU)\|_{\CS_1} &\leq C\|\bU\|_{\CS_1}^2\label{eq:LN_1}
\end{align}
for some constant $C$ depending only
 of $\mu, \nu,h$ and for any $(\bF,\bU) \in \CD_1 \times \CS_1.$

\begin{proof}[Proof of Theorem \ref{thm:U_1}]
By \cite[Proposition 2.30]{BCD2011},  we have 
\begin{equation*}
\|(-\Delta) ^{-1}\CP\bF\|_{\hB^{s_1}_{p_1,r_1}} 
\lesssim \|\bF\|_{\hB^{s_1+2}_{p_1,r_1}} 
\end{equation*}
for any $(s_1,p_1,r_1)\in \R \times [1,\infty]^2$ satisfying
$s_1+2<3\slash p_1$ or $(s_1+2,r_1)=(3\slash p_1,1).$ 
Thus \eqref{eq:L1} holds true.
\smallbreak

Next, according to the formulas in \eqref{def:N} and \cite[Proposition 2.30]{BCD2011}, it is clear that 
\begin{equation}\label{es:u_1}
\begin{aligned}
\|\CN_1 (\bU,\bU)\|_{\hB^{3/p-1}_{p,1}}
&\lesssim \mu^{-1}  \big( \|\bu \otimes \bu\|_{\hB^{3/p-2}_{p,1}}
+\|\bB \otimes \bB\|_{\hB^{3/p-2}_{p,1}} \big),\\
\|\CN_2 (\bU,\bU)\|_{\hB^{3/p-1}_{p,r}}
& \lesssim \nu^{-1} \big( \|\bu\times \bB\|_{\hB^{3/p-2}_{p,r}} +
h\|\bJ\times \bB\|_{\hB^{3/p-2}_{p,r}} \big),\\
\|\CN_3(\bU,\bU)\|_{\hB^{3/p-1}_{p,1}}
&\lesssim \nu^{-1} \big(  \|\bu\times \curl^{-1}\bJ\|_{\hB^{3/p-1}_{p,1}}
+h\| \bJ\times \curl^{-1}\bJ\|_{\hB^{3/p-1}_{p,1}} \big).
\end{aligned}
\end{equation}
Then it is easy to find from \eqref{es:pl_1} that 
\begin{equation}\label{es:u_2}
\begin{aligned}
\|\bu \otimes \bu\|_{\hB^{3/p-2}_{p,1}}
 & \lesssim \|\bu\|_{\hB^{3/p-1}_{p,2}}^2,\\
\|\bB \otimes \bB\|_{\hB^{3/p-2}_{p,1}}
&\lesssim \|\bB\|_{\hB^{3/p-1}_{p,2}}^2,\\
\| \bu \times \bB\|_{\hB^{3/p-2}_{p,r}}
&\lesssim \|\bu\|_{\hB^{3/p-1}_{p,\infty}}
\|\bB\|_{\hB^{3/p-1}_{p,r}},\\
\| \bJ \times \bB\|_{\hB^{3/p-2}_{p,r}}
&\lesssim \|\bJ\|_{\hB^{3/p-1}_{p,\infty}}
\|\bB\|_{\hB^{3/p-1}_{p,r}}
\end{aligned}
\end{equation}
provided that $1\leq p<3.$
Therefore, \eqref{es:u_1} and \eqref{es:u_2} give us the bounds of $\CN_1(\bU,\bU)$ and 
$\CN_2(\bU,\bU)$ as follows:
\begin{equation}\label{es:N1N2}
\begin{aligned}
\|\CN_1(\bU,\bU)\|_{\hB^{3\slash p-1}_{p,1}} 
&\lesssim \mu^{-1} (\|\bu\|_{\hB^{3/p-1}_{p,2}}^2
+\|\bB\|_{\hB^{3/p-1}_{p,2}}^2)
\lesssim \mu^{-1}\|\bU\|_{\CS_1}^2,\\
\|\CN_2(\bU,\bU)\|_{\hB^{3\slash p-1}_{p,r}} 
& \lesssim \nu^{-1}(\|\bu\|_{\hB^{3/p-1}_{p,\infty}}
+h\|\bJ\|_{\hB^{3/p-1}_{p,\infty}})
\|\bB\|_{\hB^{3/p-1}_{p,r}}
\lesssim \nu^{-1}(1+h)\|\bU\|_{\CS_1}^2
\end{aligned}
\end{equation}
for $r\leq 2.$
\smallbreak

On the other hand, note that 
\begin{equation*}
\|\curl^{-1}\bJ\|_{\hB^{3/p}_{p,1}} 
\lesssim \|\bJ\|_{\hB^{3/p-1}_{p,1}}.
\end{equation*}
Then \eqref{es:pl_2} implies that
\begin{align*}
\| \bu\times \curl^{-1}\bJ \|_{\hB^{3/p-1}_{p,1}} 
&\lesssim \|\bu\|_{\hB^{3/p-1}_{p,1}}\|\curl^{-1}\bJ\|_{\hB^{3/p}_{p,1}}
\lesssim \|\bu\|_{\hB^{3/p-1}_{p,1}}\|\bJ\|_{\hB^{3/p-1}_{p,1}},\\ 
 \| \bJ\times \curl^{-1}\bJ \|_{\hB^{3/p-1}_{p,1}} 
&\lesssim \|\bJ\|_{\hB^{3/p-1}_{p,1}}
\|\curl^{-1}\bJ\|_{\hB^{3/p}_{p,1}}
\lesssim \|\bJ\|_{\hB^{3/p-1}_{p,1}}^2,
\end{align*}
which together with \eqref{es:u_1} furnish that 
\begin{equation}\label{es:N3}
\|\CN_3(\bU,\bU)\|_{\hB^{3\slash p-1}_{p,1}} 
\lesssim \nu^{-1} \big(\|\bu\|_{\hB^{3/p-1}_{p,1}}
+h\|\bJ\|_{\hB^{3/p-1}_{p,1}} \big) 
\|\bJ\|_{\hB^{3/p-1}_{p,1}}
\lesssim \nu^{-1}(1+h)\|\bU\|_{\CS_1}^2.
\end{equation}
Thus \eqref{es:N1N2} and \eqref{es:N3} imply that 
\begin{equation*}
\|\CN(\bU,\bU)\|_{\CS_1} 
\leq C\big( \mu^{-1}+\nu^{-1}(1+h) \big)\|\bU\|_{\CS_1}^2
\end{equation*}
for some positive constant $C.$ This proves the estimate \eqref{eq:LN_1}.
\end{proof}

\section{Optimal well-posedness in the case  $p=2$}
This section is dedicated to the proof of Theorem \ref{thm:DS_2}, which is based on the Friedrichs' method.
The whole proof is divided into the following four steps.
\smallbreak 

\textbf{Step 1: Construction of approximate solutions.}
Set the spectral cut-off operator $\mathbb{E}_n$ by
$$\mathcal F({\mathbb{E}_n}f)(\xi)
:=\mathds{1}_{\{n^{-1}\leq|\xi|\leq n\}}(\xi)\,\mathcal F(f)(\xi).$$
Inspired by \eqref{eq:hMHD_3}, we consider the following truncated system:
\begin{equation}\label{Truncated-hMHD}
\left\{ \begin{aligned}
&  -\mu\Delta \bu= \mathbb{E}_n\mathcal P \big(\mathbb{E}_n \bB\cdot\mathbb{E}_n\nabla \bB-\mathbb{E}_n\bu\cdot\nabla \mathbb{E}_n\bu\big)+ \mathbb{E}_n\mathcal P\Bf,\\
&  -\nu\Delta \bB= \nabla\times\mathbb{E}_n 
\big( (\mathbb{E}_n\bu-h\nabla\times \mathbb{E}_n\bB)\times\mathbb{E}_n\bB \big)
+\mathbb{E}_n \CP\bg.
\end{aligned}\right.
\end{equation}
For the solvability of \eqref{Truncated-hMHD}, we recall the following lemmas:
\begin{lemma}\label{lem:fp}
Let $(X, \|\cdot\|_{X})$ be a Banach space and $\mathcal B : X\times X\to X,$  a bilinear continuous operator with norm $K$. Then, for all $a \in X$ such that $4K\|a\|_{X}<1$, the equation 
$$x=a+\mathcal B(x,x)$$ has a unique solution $x$ in the ball 
$B\big( 0, (2K)^{-1}\big)$. Besides, $x$ satisfies $\|x\|_{X}\leq2\|a\|_{X}$.
\end{lemma}

\begin{lemma} \label{lem:En}
Let $s\in\mathbb{R}, 1\leq r\leq \infty$ and  $k\geq 0$. 
Let $f\in \dot{B}_{2, r}^s(\R^3).$ Then, for all $n\geq1$, we have
\begin{align}
&\|\mathbb{E}_n f\|_{\dot{B}_{2, r}^{s+k}(\R^3)}
\lesssim n^k\|f\|_{\dot{B}_{2, r}^s(\R^3)},\label{4.100}\\
&\lim\limits_{n\to \infty}\|\mathbb{E}_n f-f\|_{\dot{B}_{2, r}^{s}(\R^3)}=0
\label{4.200},\\
&\|\mathbb{E}_n f-f\|_{\dot{B}_{2, r}^{s}(\R^3)}
\lesssim \frac{1}{n^k}\|f\|_{\dot{B}_{2, r}^{s+k}(\R^3)}\label{4.201}.
\end{align}
\end{lemma}

We claim that \eqref{Truncated-hMHD} is a system in the Banach space 
$X:=L^2(\R^3)^6$ for which Lemma \ref{lem:fp} applies.  
Indeed, introduce that 
\begin{equation*}
a:=\begin{pmatrix}
(-\mu\Delta)^{-1}\mathbb{E}_n\CP\Bf\\
 (-\nu\Delta)^{-1}\mathbb{E}_n\CP\bg
\end{pmatrix},\quad 
\CV:=\begin{pmatrix}
\bu\\
 \bB
\end{pmatrix},
\end{equation*}
\begin{equation*}
\CB(\CV,\CV):=
\begin{pmatrix}
(-\mu\Delta)^{-1}\mathbb{E}_n\mathcal P \big(\mathbb{E}_n \bB\cdot\mathbb{E}_n\nabla \bB-\mathbb{E}_n\bu\cdot\nabla \mathbb{E}_n\bu\big)\\
 (-\nu\Delta)^{-1}\nabla\times\mathbb{E}_n 
\big( (\mathbb{E}_n\bu-h\nabla\times \mathbb{E}_n\bB)\times\mathbb{E}_n\bB \big)
\end{pmatrix}.
\end{equation*}
Then Lemma \ref{lem:En} ensures that $\mathbb{E}_n$ maps $L^2(\R^3)$
to all Besov spaces, and that $\CB(\cdot,\cdot)$ is a continuous bilinear map from $X\times X$ to $X.$ 
We thus deduce that \eqref{Truncated-hMHD} admits a unique  small solution 
$(\bu^n, \bB^n)\in X,$ thanks to \eqref{small-con-p=2} and \eqref{4.100}.
Furthermore, as $\mathbb{E}_n^2=\mathbb{E}_n,$ the uniqueness implies 
$\mathbb{E}_n \bu^n=\bu^n$ and $\mathbb{E}_n \bB^n=\bB^n,$ 
and we clearly have $\dv \bu^n=\dv \bB^n=0.$
Being spectrally supported in the annulus $\{n^{-1}\leq|\xi|\leq n\},$ 
the solution $(\bu^n, \bB^n)$ belongs to $\dot B^s_{2, r}$ 
for all $s\in\R$ and $1\leq r\leq \infty.$ 
\medskip

\textbf{Step 2: Uniform estimates.} 
For $(\bu^n, \bB^n)$ obtained in last step, let us set 
$$\bv^n:=\bu^n-h\nabla\times \bB^n.$$ 
As $\dv \bv^n=\dv \bB^n=0,$ we use \cite[Identity (2.13)]{DT2021} to get
\begin{equation*}
\begin{aligned}
 \nabla\times( \bv^n \times  \bB^n )&=\bB^n \cdot\nabla  \bv^n-\bv^n\cdot\nabla \bB^n\\
 &=(\nabla \times \bv^n)\times \bB^n +(\nabla \times \bB^n)\times \bv^n 
 -2\bv^n\cdot \nabla \bB^n +\nabla (\bv^n\cdot \bB^n).
\end{aligned}
\end{equation*}
Now it is not hard to see that $(\bu^n,\bB^n,\bv^n)\in \dot{B}^{1/2}_{2, r}$ fulfils
\begin{equation}\label{Truncated-hMHD-2}
\left\{\begin{aligned}
&  -\mu\Delta\bu^n=  \mathbb{E}_n\mathcal P(\bB^n\cdot\nabla \bB^n-\bu^n\cdot\nabla \bu^n)+\mathbb{E}_n \mathcal{P}\Bf,\\
&  -\nu\Delta \bB^n= \nabla\times\mathbb{E}_n(\bv^n\times \bB^n)
+\mathbb{E}_n \CP\bg,\\
&  -\nu\Delta\bv^n
=  \mathbb{E}_n\mathcal{P}\big\{ \bB^n\cdot\nabla \bB^n-\bu^n\cdot\nabla \bu^n-h\nabla\times\big( (\nabla\times \bv^n)\times \bB^n\big)+\Bf-h\nabla\times\bg\\
&\hspace{2cm} +\nabla\times(\bv^n\times \bu^n)+2\nabla\times \big( 
\bv^n\cdot\nabla \curl^{-1}(\bu^n-\bv^n) \big)+(\mu-\nu)\Delta\bu^n\big\} \cdotp
\end{aligned}\right.
\end{equation}
The goal of this step is to derive the uniform estimates of the triplet 
$(\bu^n, \bB^n, \bv^n)$. Since $\mathbb{E}_n$ is an $L^2$ orthogonal projector, it has no effect on the energy estimates. We claim that   
\begin{equation}\label{uniform-bound-p=2}
\|(\bu^n, \bB^n, \bv^n)\|_{\dot{B}^{ {1}/{2}}_{2, r}}< 2\delta.
\end{equation}
\smallbreak

To prove \eqref{uniform-bound-p=2},  we apply the operator $\dot\Delta_j$ to both sides of \eqref{Truncated-hMHD-2}, and then we take the $L^2$ inner product with 
the vector $(\dot\Delta_j \bu^n, \dot\Delta_j \bB^n,\dot\Delta_j \bv^n)$.  
In particular, to handle the third equation of \eqref{Truncated-hMHD-2}, notice that 
$$-h\nabla\times\dot\Delta_j((\nabla\times\bv^n)\times \bB^n)\\
=\nabla\times\left([\dot\Delta_j, \curl^{-1}(\bu^n-\bv^n)\times](\nabla\times\bv^n)\right)+h\nabla\times(\bB^n\times\dot\Delta_j(\nabla\times\bv^n)),$$
and that the $L^2$ scalar product of the last term with $\dot\Delta_j\bv^n$ is 0. 
Then, we get from Bernstein lemma (see \cite[Lemma 2.1]{BCD2011} for instance) that
$$\begin{aligned}
 \mu\| \dot\Delta_j\bu^n\|_{L^2}
&\lesssim  2^{-j} \big( \|\dot\Delta_j (\bB^n \otimes \bB^n)\|_{L^2}
+\|\dot\Delta_j (\bu^n\otimes \bu^n)\|_{L^2} \big)+2^{-2j} \|\ddj \Bf\|_{L^2},\\
\nu \| \dot\Delta_j\bB^n\|_{L^2}
&\lesssim 2^{-j}\|\dot\Delta_j (\bv^n\times \bB^n)\|_{L^2}+2^{-2j}\|\ddj\bg\|_{L^2},\\
 \nu\| \dot\Delta_j\bv^n\|_{L^2}
 &\lesssim  2^{-j} \big(\|\dot\Delta_j (\bB^n \otimes \bB^n)\|_{L^2}
+\|\dot\Delta_j (\bu^n\otimes \bu^n)\|_{L^2} 
+\|\dot\Delta_j(\bv^n\times \bu^n)\|_{L^2} \big)\\
&\quad +2^{-2j} \big( \|\ddj\Bf-h\nabla\times \ddj\bg\|_{L^2}\big)
+|\mu-\nu|\,\|\ddj\bu^n\|_{L^2}\\
&\quad+2^{-j} \|\ddj(\bv^n\cdot\nabla \curl^{-1}(\bu^n-\bv^n))\|_{L^2} \\
&\quad+2^{-j} \|[\dot\Delta_j, \curl^{-1}(\bu^n-\bv^n)\times](\nabla\times\bv^n)\|_{L^2},
\end{aligned}
$$
which imply that
\begin{align*}
\mu\|\bu^n\|_{\dot{B}^{1/2}_{2, r}}
\lesssim& \|\bB^n\otimes \bB^n\|_{\dot{B}^{-1/2}_{2, r}}
+ \|\bu^n\otimes \bu^n\|_{\dot{B}^{-1/2}_{2, r}} 
+\|\Bf\|_{\dot{B}^{-3/2}_{2, r}}, \\
\nu\|\bB^n\|_{\dot{B}^{1/2}_{2, r}}
\lesssim& \|\bv^n\times \bB^n\|_{\dot{B}^{-1/2}_{2, r}}
+  \|\bg\|_{\dot{B}^{-3/2}_{2, r}}, \\
\nu\|\bv^n\|_{\dot{B}^{1/2}_{2, r}}
\lesssim& \|\bB^n\otimes \bB^n\|_{\dot{B}^{-1/2}_{2, r}}
+ \|\bu^n\otimes \bu^n\|_{\dot{B}^{-1/2}_{2, r}} 
+ \|\bv^n\times \bu^n\|_{\dot{B}^{-1/2}_{2, r}}\\
& +\|(\Bf, h\nabla \times \Bg)\|_{\dot{B}^{-3/2}_{2, r}}
+ |\mu-\nu|\|\bu^n\|_{\dot{B}^{1/2}_{2, r}}
+\|\bv^n\cdot\nabla \curl^{-1}(\bu^n-\bv^n)\|_{\dot{B}^{-1/2}_{2, r}} \\
&+\big\| 2^{-j/2}\|[\dot\Delta_j, \curl^{-1}(\bu^n-\bv^n)
\times](\nabla\times\bv^n)\|_{L^2} \big\|_{\ell^r}.
\end{align*}
According to \eqref{es:pl_1},  we see that 
\begin{equation}\label{pl:-0.5}
\|h_1 h_2\|_{\dot{B}^{-1/2}_{2, r}(\R^3)}  
\lesssim \|h_1\|_{\dot{B}^{1/2}_{2, r}(\R^3)} \|h_2\|_{\dot{B}^{1/2}_{2, r}(\R^3)}.
\end{equation}
Then by \eqref{pl:-0.5}, Proposition \ref{Prop-commu} (by taking $s=-1/2, \rho_1=\infty, \rho_2=4$)  and Young's inequality, we further obtain that 
\begin{equation*}
\|(\mu\bu^n,\nu\bB^n,\nu \bv^n)\|_{\dot{B}^{1/2}_{2, r}} 
\lesssim \bigl(1+\nu/\mu \bigr)\bigl(\|\bB^n\|_{\dot{B}^{1/2}_{2, r}}^2+ \|\bu^n\|_{\dot{B}^{1/2}_{2, r}}^2 +
 \|\bv^n\|_{\dot{B}^{1/2}_{2, r}}^2
+\|(\Bf, \bg, h\nabla\times \bg)\|_{\dot{B}^{-3/2}_{2, r}}\bigr).
\end{equation*}
Then by taking $\delta$ in \eqref{small-con-p=2} sufficiently small, we have \eqref{uniform-bound-p=2}.
\medskip
  
\textbf{Step 3: Convergence of the approximation scheme.}
To  show convergence of the nonlinear terms in system \eqref{Truncated-hMHD-2},  we need to obtain strong convergence of the solution sequence. To achieve it,  we show that the sequence  $\{(\bu^n, \bB^n, \bv^n)\}$ is a Cauchy sequence in the Banach space  $\dot{B}^0_{2, r}.$ 
  
Let $m\geq n\geq 1.$
Consider the difference $(\ddu, \ddB, \ddv):=(\bu^m-\bu^n, \bB^m-\bB^n, \bv^m-\bv^n),$  which satisfies
\begin{equation}\label{4.d1}
\left\{\begin{aligned}
 & -\mu\Delta {\ddu}=R_1+R_2,\\
 &   -\nu\Delta {\ddB}=R_3+R_4,\\
 & -\nu\Delta\ddv=R_5+R_6,
\end{aligned}
\right.
\end{equation}
with
\begin{align*}
&R_1:= (\mathbb{E}_m-\mathbb{E}_n)\mathcal P
(\bB^m\cdot\nabla \bB^m-\bu^m\cdot\nabla \bu^m+ \Bf),\\
&R_2:= \mathbb{E}_n\mathcal P(\bB^m\cdot\nabla\ddB+\ddB\cdot\nabla \bB^n-\bu^m\cdot\nabla\ddu-\ddu\cdot\nabla \bu^n),\\
&R_3:=(\mathbb{E}_m-\mathbb{E}_n) \big( \nabla\times(\bv^m\times\bB^m) \big)
+ (\mathbb{E}_m-\mathbb{E}_n)\CP\bg,\\
&R_4:=\mathbb{E}_n\nabla\times(\bv^m\times\ddB+\ddv\times \bB^n),\\
&R_5:=R_1+  (\mathbb{E}_m-\mathbb{E}_n)\Big(
\nabla\times(\bv^m\times \bu^m)-h\nabla\times \bg+(\mu-\nu)\Delta \bu^m\\
&\qquad\quad+2\nabla\times\big( \bv^m\cdot\nabla\curl^{-1}(\bu^m-\bv^m) \big) 
-h\nabla\times \big( (\nabla\times\bv^m)\times \bB^m\big)\Big),\\
&R_6:=R_2+ \mathbb{E}_n\nabla\times\bigl( \bv^m\times\ddu+\ddv\times \bu^n-h(\nabla\times\bv^m)\times\ddB-h(\nabla\times\ddv)\times \bB^n\\
&\qquad\quad +2\bv^m\cdot\nabla\curl^{-1}(\ddu-\ddv)
+2\ddv\cdot\nabla\curl^{-1}( \bu^n-\bv^n)\bigr)+(\mu-\nu)\Delta \ddu.
\end{align*}
Hence, arguing as in the second step of the proof and using the following decomposition
$$\begin{aligned}
&-h\ddj\mathbb{E}_n\nabla\times ((\nabla\times\ddv)\times \bB^n)\\
=&\mathbb{E}_n\nabla\times \big( [\ddj, \curl^{-1}(\bu^n-\bv^n) \times](\nabla\times \ddv) \big)+ h\mathbb{E}_n\nabla\times \big( \bB^n\times(\nabla\times \ddj\ddv) \big)
\end{aligned}$$
give us  that
\begin{equation}\label{es:d1}
\|(\mu\ddu, \nu\ddB, \nu\ddv)\|_{\dot{B}^{0}_{2, r}}
\lesssim (1+\nu/\mu) \Bigl(\|(R_1, R_2,  R_3, R_4, R_5)\|_{\dot{B}^{-2}_{2, r}} 
+\sum_{j=1}^3 K_j \Bigr),
\end{equation}
where we have set 
\begin{align*}
K_1 &:=\| \bv^m\times \ddu \|_{\dot{B}^{-1}_{2, r}}
+\|\ddv\times \bu^n \|_{\dot{B}^{-1}_{2, r}}\\
&\quad +\|\bv^m\cdot\nabla\curl^{-1}(\ddu-\ddv)\|_{\dot{B}^{-1}_{2, r}}
+\|\ddv\cdot\nabla\curl^{-1}( \bu^n-\bv^n)\|_{\dot{B}^{-1}_{2, r}},\\
K_2 &:= h\| (\nabla\times \bv^m)\times\ddB\|_{\dot{B}^{-1}_{2, r}},\\
K_3 &:=\big\|\|2^{-j}[\ddj, \curl^{-1}(\bu^n-\bv^n)\times] 
(\nabla\times\ddv)\|_{L^2} \big\|_{\ell^r}.
\end{align*}

At this stage, most terms in $R_1, R_3, R_5$ can be estimated by \eqref{4.201} and \eqref{pl:-0.5} except for the last term of $R_5.$ 
Indeed, we have 
\begin{equation}\label{es:R5_1}
\begin{aligned}
\|(R_1, R_3, R_5)\|_{\dot{B}^{-2}_{2, r}}
& \lesssim \frac{1}{\sqrt{n}}\Bigl(
\|\bB^m\|_{\dot{B}^{1/2}_{2, r}}^2+\|\bu^m\|_{\dot{B}^{1/2}_{2, r}}^2
+\|\bv^m\|_{\dot{B}^{1/2}_{2, r}}^2
+\|(\Bf, \bg, h\nabla\times\bg)\|_{\dot{B}^{-3/2}_{2, r}}\Bigr)  \\
& \quad  +\frac{|\mu-\nu|}{\sqrt{n}}\,\|\bu^m\|_{\dot{B}^{1/2}_{2, r}}
+h \big\| (\mathbb{E}_m-\mathbb{E}_n)  \big( (\nabla\times\bv^m)\times \bB^m\big)\big\|_{\dot{B}^{-1}_{2, r}}.
\end{aligned}
\end{equation}
According to the continuity of paraproduct and remainder operator in Besov spaces, 
the following product law holds 
\begin{equation}\label{pl:-0.75}
\|h_1 h_2\|_{\dot{B}^{-3/4}_{2, r}(\R^3)}  
\lesssim \|h_1\|_{\dot{B}^{-1/2}_{2, r}(\R^3)} \|h_2\|_{\dot{B}^{5/4}_{2, r}(\R^3)}.
\end{equation}
Then we use \eqref{4.201}, \eqref{pl:-0.75} and the interpolation theory to obtain 
\begin{align*}
\big\| (\mathbb{E}_m-\mathbb{E}_n)  
\big( (\nabla\times\bv^m)\times \bB^m\big)\big\|_{\dot{B}^{-1}_{2, r}} 
 \lesssim &  \frac{1}{n^{1/4}} \| (\nabla\times\bv^m)\times \bB^m\|_{\dot{B}^{-3/4}_{2, r}} \\
 \lesssim & \frac{1}{n^{1/4}}\|\bv^m\|_{\dot{B}^{1/2}_{2, r}}\|\bB^n\|_{\dot{B}^{5/4}_{2, r}} \\
 \lesssim & \frac{1}{n^{1/4}}\|\bv^m\|_{\dot{B}^{1/2}_{2, r}}\|\bB^n\|_{\dot{B}^{1/2}_{2, r}}^{1/4}\big\|\frac{1}{h}(\bu^n-\bv^n) \big\|_{\dot{B}^{1/2}_{2, r}}^{3/4}\\
 \lesssim &\frac{1}{n^{1/4} h^{3/4}} \|\bv^m\|_{\dot{B}^{1/2}_{2, r}}\|\bB^n\|_{\dot{B}^{1/2}_{2, r}}^{1/4}\| (\bu^n, \bv^n)\|_{\dot{B}^{1/2}_{2, r}}^{3/4}.
\end{align*}
Inserting the above bound into \eqref{es:R5_1}, we obtain that 
\begin{equation}\label{es:R135}
\begin{aligned}
\|(R_1, R_3, R_5)\|_{\dot{B}^{-2}_{2, r}}
& \lesssim \frac{1}{\sqrt{n}}\Bigl(
\|\bB^m\|_{\dot{B}^{1/2}_{2, r}}^2+\|\bu^m\|_{\dot{B}^{1/2}_{2, r}}^2
+\|\bv^m\|_{\dot{B}^{1/2}_{2, r}}^2
+\|(\Bf, \bg, h\nabla\times\bg)\|_{\dot{B}^{-3/2}_{2, r}}\Bigr)  \\
& \quad  +\frac{|\mu-\nu|}{\sqrt{n}}\,\|\bu^m\|_{\dot{B}^{1/2}_{2, r}}
+\Big( \frac{h}{n} \Big)^{1/4}  \|\bv^m\|_{\dot{B}^{1/2}_{2, r}}\|\bB^n\|_{\dot{B}^{1/2}_{2, r}}^{1/4}\| (\bu^n, \bv^n)\|_{\dot{B}^{1/2}_{2, r}}^{3/4}.
\end{aligned}
\end{equation}
\smallbreak 
 
Next, we shall prove that 
\begin{equation}\label{es:R24K}
\|(R_2,R_4)\|_{\dot{B}^{-2}_{2, r}} +\sum_{j=1}^3 K_j 
\lesssim  \| (\bu^m, \bu^n,\bB^m, \bB^n, \bv^m,\bv^n)\|_{\dot{B}^{1/2}_{2, r}}   
\|(\ddu,\ddB,\ddv)\|_{\dot{B}^{0}_{2, r}}.
\end{equation}
Using the following product law
\begin{equation*}
\|h_1 h_2\|_{\dot{B}^{-1}_{2, r}(\R^3)}  
\lesssim \|h_1\|_{\dot{B}^{0}_{2, r}(\R^3)} \|h_2\|_{\dot{B}^{1/2}_{2, r}(\R^3)},
\end{equation*}
it is not hard to see that 
\begin{equation}\label{es:R24}
\begin{aligned}
&\|(R_2,R_4)\|_{\dot{B}^{-2}_{2, r}} +K_1 
\lesssim  \| (\bu^m, \bu^n,\bB^m, \bB^n, \bv^m,\bv^n)\|_{\dot{B}^{1/2}_{2, r}}   
\|(\ddu,\ddB,\ddv)\|_{\dot{B}^{0}_{2, r}}.
\end{aligned}
\end{equation}
On the other hand, the bound 
\begin{equation*}
\|h_1 h_2\|_{\dot{B}^{-1}_{2, r}(\R^3)}  
\lesssim \|h_1\|_{\dot{B}^{-1/2}_{2, r}(\R^3)} \|h_2\|_{\dot{B}^{1}_{2, r}(\R^3)}
\end{equation*}
implies that   
\begin{equation}\label{es:K2}
K_2 \lesssim  h\|\nabla\times \bv^m\|_{\dot{B}^{-1/2}_{2, r}} \|\ddB\|_{\dot{B}^{1}_{2, r}}\lesssim \| \bv^m\|_{\dot{B}^{1/2}_{2, r}} \|(\ddu, \ddv)\|_{\dot{B}^{0}_{2, r}}.
\end{equation}
To study $K_3$, we use Proposition \ref{Prop-commu} (by taking $s=-1, \rho_1=\infty, \rho_2=2$) to see
\begin{equation}\label{es:K3}
 \begin{aligned}
K_3\lesssim & \| \curl^{-1}(\bu^n-\bv^n)\|_{\dot{B}^{3/2}_{2, \infty}}\|\nabla\times\ddv\|_{\dot{B}^{-1}_{2, r}} +   \| \curl^{-1}(\bu^n-\bv^n)\|_{\dot{B}^{3/2}_{2, r}}\|\nabla\times\ddv\|_{\dot{B}^{-1}_{2, \infty}} \\
\lesssim & \|(\bu^n, \bv^n)\|_{\dot{B}^{1/2}_{2, r}}\|\ddv\|_{\dot{B}^{0}_{2, r}}.
 \end{aligned}
 \end{equation}
Thus \eqref{es:R24}, \eqref{es:K2} and \eqref{es:K3} furnish \eqref{es:R24K}.
\medskip
 
Now,  keeping \eqref{small-con-p=2} and \eqref{uniform-bound-p=2} in mind,  we insert \eqref{es:R135} and \eqref{es:R24K} into \eqref{es:d1} to get 
\begin{equation}\label{es:d2}
\|(\ddu, \ddB, \ddv)\|_{\dot{B}^{0}_{2, r}}
\lesssim n^{-1/4}
\end{equation}
from which we conclude that $\{(\bu^n, \bB^n, \bv^n)\}$ is a Cauchy sequence in $\dot{B}^{0}_{2, r}(\R^3).$
\medskip

\textbf{Step 4: Existence and uniqueness.}
As a consequence of  \eqref{uniform-bound-p=2}, \eqref{es:d2} and  the interpolation theory, $(\bu^n, \bB^n)$ converges strongly to some  $(\bu, \bB)$  in all $\dot{B}^{s_1}_{2, r}\times  \dot{B}^{s_2}_{2, r}$ for $0\leq s_1< 1/2$ and $0\leq s_2< 3/2.$ This shows that the limit $(\bu, \bB)$ solves original system \eqref{eq:hMHD_1} in the sense of distributions. Uniqueness is guaranteed by performing similar estimates like in the third step. The details are omitted.
This completes the proof of Theorem \ref{thm:DS_2}.

\section{Ill-posedness issue in critical Besov space}
\label{sec:ill}
In this section, we prove Theorem \ref{thm:DS_3}.
Let us consider the stationary Naver-Stokes equations \eqref{eq:NS_1} in the following form:
\begin{equation}\label{eq:NS_2}
\bu= L\Bf + N(\bu,\bu),
\end{equation}
where 
$$L\Bf = (-\mu\Delta)^{-1}\CP \Bf
\,\,\,\text{and}\,\,\,
N(\bu,\bu)=-(-\mu\Delta)^{-1} \CP \dv( \bu \otimes \bu).$$

\subsection{Case $p>3$ and $p=3$ with $r>2$}
In this subsection, we shall review that the solution mapping of \eqref{eq:NS_2} is discontinuous from $\big(B_{\wt{\CD}}(\ep),\|\cdot\|_{\wt{\CD}}\big)$ to $\big(\wt{\CS}, \|\cdot\|_{\wt{\CS}}\big)$ with
\begin{equation}\label{def:tilde}
\wt{\CD} :=\hB^{3/\wt p-3}_{\wt p,r}(\R^3), \quad (\wt p, r)\in ((3,\infty]\times[1,\infty])\cup(\{3\}\times(2,\infty]),\quad 
\wt{\CS} :=\hB^{-1}_{\infty, \infty}(\R^3).
\end{equation}
In fact, \cite{Ts2019} proved the following results.
\begin{lem}\label{lem:ill}
Let $\wt \CD$ and $\wt \CS$ be given in \eqref{def:tilde}.
Moreover, we define the subspace $\CD\subset \wt \CD$ as
\begin{equation*}
\CD:=\begin{cases}
\hB^{-2}_{3,1}(\R^3) & {\rm if}\,\,\, \wt p>3,\\
\hB^{-2}_{3,2}(\R^3) &  {\rm if}\,\,\, \wt p=3.
\end{cases}
\end{equation*}
Then for any small $\ep>0,$ there exist a sequence 
$\{\Bf_n\}_{n\in\mathbb{N}}$
 of external forces and a constant $c_0=c_0(\ep)>0$ satisfying 
\begin{enumerate}
\item $\|\Bf_n\|_{\CD}<\ep$\quad for every $n\in\mathbb{N}$;
\item $\|\Bf_n\|_{\wt \CD}\to 0$ as $n\to\infty$;
\item $\|N(L\Bf_n,L\Bf_n)\|_{\hB^{-1}_{\infty,\infty}}>c_0,$ so long as $n$ large enough.
\end{enumerate}
\end{lem}

\begin{rmk}
Suppose that Lemma \ref{lem:ill} holds here.
Then clearly, we have $\|\Bf_n\|_{\wt\CD}<\ep$ for every $n\in\mathbb{N}$.
Moreover, since \eqref{eq:NS_2} is well-posed from $\CD$ to $L^3(\R^3)$ by \cite{Ts2019},
\eqref{eq:NS_2} admits a solution $\bu_n\in L^3(\R^3)\subset \wt\CS$ for such force $\Bf_n$ $(n\in\mathbb{N})$ with taking sufficiently small $\ep$.
However, Lemma \ref{lem:ill} and the contraposition of Proposition \ref{prop:BT} (by the lack of continuity of $A_2$) imply that the convergence $\Bf_n\to 0$ in a weaker space $\wt\CD$ does not necessarily guarantee the convergence $\bu_n \to 0$ even in the weakest space $\wt\CS$.
This means the ill-posedness of the system \eqref{eq:NS_2} (and the system \eqref{eq:hMHD_1} with $\bg\equiv \bB\equiv 0$) from $\wt\CD$ to $\wt\CS$.
\end{rmk}

\subsection{Case $p=3$ with $r\leq2$}
In this subsection, we will prove some results on the well-posedness 
and ill-posedness of  three-dimensional stationary Navier-Stokes equations in the critical Besov spaces for the index $p=3.$
Let us introduce 
\begin{equation*}
\CD_3:=\hB^{-2}_{3,r}(\R^3)^3,
\quad
\CS_3:=\big\{\bu\in \hB^{0}_{3,r} (\R^3)^3:\dv \bu=0\big\}.
\end{equation*}
Then the goal of this subsection is to prove the following results.
\begin{thm}\label{thm:NS_1}
Let $3/2<r\leq 2.$ 
Then the system \eqref{eq:NS_2} is well-posed from $\CD_{3}$ to $\CS_{3}.$ 
\end{thm}
\begin{thm}\label{thm:NS_3}
Let $1\leq r<3/2.$ Then the system \eqref{eq:NS_2} is ill-posed from $\CD_{3}$ to $\CS_{3}.$ 
\end{thm}

\subsubsection{Proof of Theorem \ref{thm:NS_1}}
To prove Theorem \ref{thm:NS_1}, we recall the definition and the property of the homogeneous Triebel-Lizorkin space  $\hF^s_{p,r}(\R^3).$
\begin{dfn} 
Let $s\in\R,$ $(p, r)\in [1,\infty)\times[1,\infty]$ or $p=r=\infty.$
 $\hF^s_{p,r}(\R^3)$ is defined by
\begin{equation*}\label{def:TL}
\hF^{s}_{p,r}(\R^3):= \big\{ \,u \in  \mathcal{S}'(\R^3)/\CP(\R^3): 
\|u\|_{\hF^s_{p,r}} :=\big\|\|2^{js}\hDelta_j u \|_{\ell^r(\mathbb{Z}) }
 \big\|_{L^p(\R^3)} <\infty \,\big\}.
\end{equation*}
In addition, if $1<p<\infty,$ then $\hF^0_{p, 2} (\R^3)$ coincides with the Lebesgue space, that is, 
\begin{equation}\label{eq:withL}
\hF^0_{p, 2} (\R^3)\simeq L^p(\R^3).
\end{equation}
\end{dfn}

\begin{rmk}
The definition of the homogeneous Triebel-Lizorkin space  $\hF^s_{\infty,r}(\R^3)$ ($1\leq r<\infty$) is more involved, and the reader may refer to  \cite[Chapter 2]{Gr2014} and \cite[Chapter 2]{Sa2018} for details.
\end{rmk}

\begin{prop} 
Let $s\in \R$ and $1\leq r_1 \leq r_2\leq \infty.$  The following embeddings hold true.
\begin{enumerate}
\item For any $1\leq p \leq \infty$ and $1\leq p_1<p_2 \leq \infty,$ we have 
\begin{align}\label{eq:embd_TL}
\hF^{s}_{p,r_1}(\R^3) \hookrightarrow 
\hF^{s}_{p,r_2}(\R^3), \quad 
\hF^{s}_{p_1,\infty}(\R^3)  \hookrightarrow 
\hF^{s-3(1\slash p_1-1\slash p_2)}_{p_2,1}(\R^3) .
\end{align}

\item For any $(p, r)\in [1,\infty]^2,$ we  have 
\begin{equation}
\label{eq:withB}
\hB^s_{p, \min(p,r)} (\R^3)\hookrightarrow
\hF^s_{p, r} (\R^3)\hookrightarrow
\hB^s_{p, \max(p,r)}(\R^3).
\end{equation}
\end{enumerate}

\end{prop}

Theorem \ref{thm:NS_1} is immediately from the following lemma and Proposition \ref{prop:BT}.
\begin{lem}\label{lem:NNS_1}
Let $3/2<r\leq 2.$ Then we have 
\begin{equation*}
\|N(\bu,\bv)\|_{\hB^{0}_{3,r}(\R^3)} 
\lesssim \|\bu\|_{L^{3}(\R^3)}  \|\bv\|_{L^{3}(\R^3)} 
\lesssim  \|\bu\|_{\hB^{0}_{3,r}(\R^3)} \|\bv\|_{\hB^{0}_{3,r}(\R^3)}.
\end{equation*}
\end{lem}
\begin{proof}
Using \eqref{eq:withL}, \eqref{eq:embd_TL} and \eqref{eq:withB}, we easily have
\begin{equation*}
\begin{aligned}
\|N(\bu,\bv)\|_{\hB^{0}_{3,r}} \lesssim 
\|\bu \otimes \bu\|_{\hB^{-1}_{3,r}}
&\lesssim \|\bu \otimes \bu\|_{\hB^{3/r-2}_{r,r}}
\lesssim \|\bu \otimes \bu\|_{\hF^{3/r-2}_{r,r}}\\
&\lesssim \|\bu \otimes \bu\|_{\hF^{0}_{3/2,2}}
\lesssim \|\bu \otimes \bu\|_{L^{3/2}}\\
&\lesssim \|\bu\|_{L^{3}}^2
\lesssim \|\bu\|_{\hF^{0}_{3,2}}^2
\lesssim \|\bu\|_{\hB^{0}_{3,2}}^2
\lesssim \|\bu\|_{\hB^{0}_{3,r}}^2.
\end{aligned}
\end{equation*}
\end{proof}

In addition, for the case where $1\leq r\leq 3/2,$ we have the following technical result which will be applied to study the ill-posedness issue.
\begin{thm}\label{thm:NS_2}
Let $0<\sigma_1<1,$ $0<\sigma_2\leq 1/2$  and $1\leq r \leq 3/2.$  
Suppose that 
\begin{equation*}
p_{1}=3/(1+\sigma_1),\quad 
r_2=3/(2-\sigma_2)
\quad \text{and} \quad 
1/r=1/r_1 + 1/r_2.
\end{equation*}
Let us define
\begin{equation*}
\CD_{\sigma} :=\hB^{-2}_{3,r}(\R^3)^3 \cap \hB^{\sigma_1-2}_{p_{1},r_1}(\R^3)^3
\,\,\,\text{and}\,\,\,
\CS_{\sigma} :=\hB^{0}_{3,r}(\R^3)^3 \cap \hB^{\sigma_1}_{p_{1},r_1}(\R^3)^3.
\end{equation*}
Suppose that $\Bf\in \CD_{\sigma}$ with 
$$ \|\Bf\|_{\hB^{-2}_{3,r_2}} \leq \delta $$
for some small constant $\delta.$ Let $\bu$ be the solution given by Theorem \ref{thm:NS_1} satisfying 
\begin{equation*}
\|\bu\|_{\hB^{0}_{3,r_2}} \leq C \delta.
\end{equation*}
Then we have $\bu \in \CS_{\sigma}.$
\end{thm}

\begin{proof}
The estimates of $L\Bf$ are clear, that is,
\begin{equation*}
\|L\Bf\|_{\hB^{0}_{3,r}} 
\lesssim \|\Bf\|_{\hB^{-2}_{3,r}},\quad 
\|L\Bf\|_{\hB^{\sigma_1}_{p_{1},r_1}}
 \lesssim \|\Bf\|_{\hB^{\sigma_1-2}_{p_{1},r_1}}.
\end{equation*}
Let $q_{1}:=3/(2+\sigma_1)<p_{1}$ and $r_2=3/(2-\sigma_2).$ 
 In particular, we have $r_1<r.$
By the Bony Calculus, we have 
\begin{equation*}
\begin{aligned}
\|\dot T_u  v \|_{\hB^{-1}_{3,r}} 
&\lesssim  \|u\|_{\hB^{-1}_{\infty,r}} \|v\|_{\hB^{0}_{3,\infty}}
\lesssim  \|u\|_{\hB^{0}_{3,r}} \|v\|_{\hB^{0}_{3,r_2}}, \\
\|\dot T_v  u \|_{\hB^{-1}_{3,r}} 
&\lesssim  \|v\|_{\hB^{-1}_{\infty,\infty}} \|u\|_{\hB^{0}_{3,r}}
\lesssim  \|v\|_{\hB^{0}_{3,r_2}} \|u\|_{\hB^{0}_{3,r}}, \\
\|\dot R(u,v)\|_{\hB^{-1}_{3,r}} 
&\lesssim \|\dot R(u,v)\|_{\hB^{\sigma_1}_{q_{1},r}} 
\lesssim \|u\|_{\hB^{\sigma_1}_{p_{1},r_1}} \|v\|_{\hB^{0}_{3,r_2}},\\
\|\dot T_{u} v \|_{\hB^{\sigma_1-1}_{p_{1},r_1}} 
&\lesssim  \|u\|_{\hB^{\sigma_1-1}_{3/\sigma_1,r_1}} \|v\|_{\hB^{0}_{3,\infty}}
\lesssim \|u\|_{\hB^{\sigma_1}_{p_{1},r_1}}\|v\|_{\hB^{0}_{3,r_2}}, \\
\|\dot T_{v} u \|_{\hB^{\sigma_1-1}_{p_{1},r_1}} 
&\lesssim \|v\|_{\hB^{-1}_{\infty,\infty}}  \|u\|_{\hB^{\sigma_1}_{p_1,r_1}}
\lesssim \|v\|_{\hB^{0}_{3,r_2}}  \|u\|_{\hB^{\sigma_1}_{p_1,r_1}} , \\
\|\dot R(u,v)\|_{\hB^{\sigma_1-1}_{p_{1},r_1}} 
&\lesssim \|\dot R(u,v)\|_{\hB^{\sigma_1}_{q_{1},r_1}} 
\lesssim  \|u\|_{\hB^{\sigma_1}_{p_1,r_1}}
\|v\|_{\hB^{0}_{3,\infty}} , 
\end{aligned}
\end{equation*}
which imply that 
\begin{equation*}
\begin{aligned}
 \|N(\bu,\bv)\|_{\hB^{0}_{3,r}} &\lesssim \|\bu \otimes \bv\|_{\hB^{-1}_{3,r}}
 \lesssim \|\bv\|_{\hB^{0}_{3,r_2}} \big( \|\bu\|_{\hB^{0}_{3,r}}
 + \|\bu\|_{\hB^{\sigma_1}_{p_{1},r_1}} \big),\\
  \|N(\bu,\bv)\|_{\hB^{\sigma_1}_{p_{1},r_1}} 
 & \lesssim \|\bu \otimes \bv\|_{\hB^{\sigma_1-1}_{p_{1},r_1}}
 \lesssim \|\bv\|_{\hB^{0}_{3,r_2}} \|\bu\|_{\hB^{\sigma_1}_{p_{1},r_1}} .
\end{aligned}
\end{equation*}
Thus \eqref{eq:NS_2} yields that
\begin{equation*}
\begin{aligned}
\|\bu\|_{\CS_{\sigma}} 
\leq C \|\Bf\|_{\CD_{\sigma}} 
+  \|N(\bu,\bu)\|_{\CS_{\sigma}} 
\leq C  \|\Bf\|_{\CD_{\sigma}}  + C \delta\|\bu\|_{\CS_{\sigma}}.
\end{aligned}
\end{equation*}
As $\delta$ is small, we have 
\begin{equation*}
\|\bu\|_{\CS_{\sigma}} \leq C \|\Bf\|_{\CD_{\sigma}}.
\end{equation*}
\end{proof} 

\subsubsection{Proof of Theorem \ref{thm:NS_3}}

Suppose that
\begin{equation}\label{eq:r}
1\leq r=3/(2+\ep)<3/2
\end{equation}
for some $0<\ep\leq 1.$ 
We shall prove the norm inflation phenomenon of \eqref{eq:NS_2} in the space $\hB^{0}_{3,r}(\BR^3)$ for such $r$ by using the forces of \cite{LYZ2022} in three-dimensional setting. 
More precisely, let us set
\begin{equation*}
16 \BN= \{16 k: k=1,2,\cdots\}, \quad \BN(n)=\{ k \in 8\BN: n/4 \leq k \leq n/2\}.
\end{equation*}
For such set $\BN(n),$ we set the semi-norm
\begin{equation*}
\|u\|_{\hB^{0}_{3,q} (\BN(n))} := \Big( \sum_{j \in\BN(n) }\|\hDelta_j u \|_{L^3(\R^3)}^q   \Big)^{1/q} \quad \text{for any}\,\,\,1\leq q<\infty.
\end{equation*}
\smallbreak 

In addition, let $\hat \theta \in C_0^{\infty}(\BR)$ satisfy 
\begin{equation*}
0\leq \wh \theta (\cdot) \leq 1
\quad \text{and}\quad 
\wh  \theta (\xi) = 
\begin{cases}
1 & \text{if}\,\,\, |\xi|\leq 1/600,\\
 0 & \text{if}\,\,\,  |\xi| \geq 1/300.
\end{cases}
\end{equation*}
Then we define
\begin{equation*}
\phi (x)=\theta(x_1) \theta(x_2)\theta(x_3)
\sin \Big(\frac{17}{24} x_3\Big).
\end{equation*}
Moreover, for any $n\in \BN,$ we define
\begin{equation}\label{def:bc}
\begin{aligned}
b_n &:= \,n^{-1/(2r)} \sum_{k \in \BN(n)} 2^{k} \phi \big( 2^k A(x-2^{2n+k} \be) \big) 
\sin \Big(\frac{17}{12} 2^n (x\cdot \be)\Big)\\
c_n&:= \CF^{-1}\Big[ \frac{(\xi_2-\xi_1)\wh b_n(\xi)}{\xi_2}\Big],\qquad 
\bg_n=\begin{pmatrix}
b_n\\
c_n-b_n\\
0
\end{pmatrix}.
\end{aligned}
\end{equation}
with 
\begin{equation*}
\be:=\frac{\sqrt{2}}{2} (1,1,0)^{\top}, \quad 
\text{and}\quad 
A= \begin{pmatrix}
 \ep  & 0 & 0\\
0 &   \ep & 0\\
0 & 0 &  1\\
\end{pmatrix}.
\end{equation*}

Furthermore, we claim that the bounds of $b_n,$ $c_n$ and $\bg_n$ obtained in \cite{LYZ2022} still hold true for three-dimensional case. More precisely, we have the following technical results. 
\begin{lem}\label{lem:key}
Let $r$ be given by \eqref{eq:r}, and let $b_n,$ $c_n$ and $\bg_n$ be defined by \eqref{def:bc}.
The following assertions hold true so long as $n$ large enough.
\begin{enumerate}[$(1)$]

\item For $\fh_n\in \{b_n,c_n,\bg_n\},$ we have 
\begin{equation*}
\hDelta_j \fh_n=\begin{cases}
\fh_n & \text{for}\,\,\,j=n,\\
0 & \text{for}\,\,\,j\not=n.
\end{cases}
\end{equation*}
Moreover, the following estimates hold
\begin{equation*}
\begin{aligned}
\|b_n\|_{L^3 (\BR^3)} +\|b_n\|_{\hB^{0}_{3,1} (\BR^3)}
&\lesssim n^{-\ep/6},\\
\|c_n\|_{L^3(\BR^3)}+\|c_n\|_{\hB^{0}_{3,1} (\BR^3)}
&\lesssim 2^{-n/3},\\
\|\bg_n\|_{L^3 (\BR^3)} +\|\bg_n\|_{\hB^{0}_{3,1} (\BR^3)}
&\lesssim n^{-\ep/6}.
\end{aligned}
\end{equation*}

\item Setting that
\begin{equation}\label{def:G}
\bG_n:=N(\bg_n, \bg_n)
=-(-\mu\Delta)^{-1} \CP \dv( \bg_n \otimes \bg_n),
\end{equation}
we have 
\begin{equation*}
 \|\bG_n\|_{\hB^{0}_{3,r} (\BN(n))} \geq c
\end{equation*}
for some positive constant $c.$
\end{enumerate}
\end{lem}
Here, we omit the details of the proof of Lemma \ref{lem:key} for simplicity.
Using Lemma \ref{lem:key}, it is easy to see the following result.
\begin{lem} \label{lem:bc}
Let $1\leq p\leq 3$ and let $r$ be given by \eqref{eq:r}.  Then we have 
\begin{equation*}
\|b_n\|_{L^p(\BR^3)} \lesssim n^{1/p-1/(2r)} 2^{-n(3/p-1)/2},\qquad
\|c_n\|_{L^p(\BR^3)} \lesssim 2^{-n(11/p-1)/8}.
\end{equation*}
\end{lem}
\begin{proof}
As $\phi$ is a Schwarz function, we have 
\begin{equation*}
|\phi(x)|+|\nabla \phi (x)| \lesssim (1+|x|)^{-4}.
\end{equation*}
Then it is easy to see that 
\begin{equation*}
\begin{aligned}
  \|b_n\|_{L^1(\BR^3)} 
\leq &n^{-1/2r}\sum_{k\in \BN(n)}\int_{\BR^3} 
\frac{2^{k} \,dx}{(1+|2^{k} A(x-2^{2n+k}\be)|)^4} \\
\leq &n^{-1/2r}\sum_{k\in \BN(n)}\int_{\BR^3} 
\frac{2^{k} \,dx}{(1+|2^{k} Ax|)^4} \\
\leq &n^{-1/2r}\sum_{k\in \BN(n)}2^{-2k} \ep^{-2}\int_{\BR^3} 
\frac{dx}{(1+|x|)^4} \\
\lesssim & \ep^{-2}  n^{1-1/2r} 2^{-n}.
\end{aligned}
\end{equation*}
Hence the interpolation theory and Lemma \ref{lem:key} yield that
\begin{equation*}
\|b_n\|_{L^p (\BR^3)}  \lesssim 
\|b_n\|_{L^1(\BR^3)}^{(3/p-1)/2} 
\|b_n\|_{L^3(\BR^3)}^{(3-3/p)/2}
\lesssim n^{1/p-1/(2r)} 2^{-n(3/p-1)/2}.
\end{equation*}

Now, let us set
\begin{equation*}
\Phi_k(\xi)
:=\phi \big( 2^k A(x-2^{2n+k} \be) \big) 
\sin \Big(\frac{17}{12} 2^n (x\cdot \be)\Big).
\end{equation*}
Near the support of $\wh \Phi_k(\xi),$ according to the choice of $\theta$, we can prove 
\begin{equation*}
|\xi_2-\xi_1|\lesssim \ep 2^k, \quad 2^n \lesssim |\xi_2|.
\end{equation*}
Then we have 
\begin{equation*}
\begin{aligned}
\|c_n\|_{L^1(\R^3)} 
& \lesssim \|(\xi_2-\xi_1)\wh b_n(\xi)/\xi_2 \|_{L^{\infty}(\R^3)} \\
& \lesssim \ep  2^{-n}
\sup_{k\in \BN(n)} \Big( 2^{2k}\big\| \wh \Phi_k(\xi) \big\|_{L^{\infty}(\R^3)}  \Big)\\
& \lesssim \ep 2^{-n} \sup_{k\in \BN(n)} (\ep^{-2} 2^{-k})\\
& \lesssim  \ep^{-1} 2^{-5n/4}.
\end{aligned}
\end{equation*}
Thus, by the interpolation theory and Lemma \ref{lem:key}, we obtain that 
\begin{equation*}
\|c_n\|_{L^p (\BR^3)}  \lesssim 
\|c_n\|_{L^1(\BR^3)}^{(3/p-1)/2} 
\|c_n\|_{L^3(\BR^3)}^{(3-3/p)/2}
\lesssim 2^{-n(11/p-1)/8}.
\end{equation*}
\end{proof}

For  $r$ given by \eqref{eq:r},
we choose the indices $\sigma_1,$ $p_1$ $r_1$ and $r_2$ as in Theorem \ref{thm:NS_2} 
additionally satisfying the condition
\begin{equation}\label{cdt:ss}
0<2\max\{\sigma_1, \sigma_2\}<\ep\leq 1.
\end{equation}
If we take the forces 
\footnote{Using the equality $(\pd_1-\pd_2)b_n+\pd_2 c_n=0,$
then we immediately have $\dv \bg_n=0,$ and thus $\dv \Bf_n=0.$}
$\Bf_n:=-\mu\Delta\bg_n,$ then Lemmas \ref{lem:key} and \ref{lem:bc} imply that 
\begin{equation}\label{es:force}
\begin{aligned}
\|\bg_n\|_{L^{p_{1}}(\R^3)} 
&\lesssim n^{(2\sigma_1-\ep)/6}2^{-n\sigma_1/2}+2^{-n(11\sigma_1+8)/24}
\lesssim n^{(2\sigma_1-\ep)/6}2^{-n\sigma_1/2},\\
\|\Bf_n\|_{\hB^{-2}_{3,1}(\R^3)}
&\lesssim \|\bg_n\|_{\hB^{0}_{3,1}(\R^3)}
\lesssim n^{-\ep/6},\\
\|\Bf_n\|_{\hB^{\sigma_1-2}_{p_{1},r_1}(\R^3)} 
&\lesssim \|\bg_n\|_{\hB^{\sigma_1}_{p_{1},r_1}(\R^3)}
\lesssim 2^{n\sigma_1} \|\bg_n\|_{L^{p_{1}}(\R^3)}
\lesssim n^{(2\sigma_1-\ep)/6}2^{n\sigma_1/2}.
\end{aligned}
\end{equation}
Here, we also used the condition \eqref{cdt:ss}.
Thanks to Theorem \ref{thm:NS_2}, the system \eqref{eq:NS_2} admits a solution 
$$\bu_n \in \CS_{\sigma}=\hB^{0}_{3,r}(\R^3)^3 \cap \hB^{\sigma_1}_{p_{1},r_1}(\R^3)^3$$ associated with such forces $\Bf_n.$
Moreover, we have
\begin{equation*}
\|\bu_n\|_{\hB^{0}_{3,r_2}(\R^3)} \lesssim n^{-\ep/6} 
\quad \text{for}\,\,\, r_2=3/(2-\sigma_2).
\end{equation*}
Denote 
$$\bU_n:=\bu_n-\bg_n-\bG_n.$$
with $\bG_n$ defined in \eqref{def:G}.
From \eqref{eq:NS_2}, we have
\begin{equation*}
\begin{aligned}
\bU_n = & N(\bU_n, \bg_n)+N(\bg_n,\bU_n)
+N(\bU_n,\bG_n)+N(\bU_n,\bG_n)\\
&+N(\bG_n, \bg_n)+N(\bg_n,\bG_n)
+ N(\bU_n,\bU_n)+N(\bG_n,\bG_n).
\end{aligned}
\end{equation*}
Using the above equation of $\bU_n$ and \eqref{es:force}, we can prove that 
\begin{equation*}
\begin{aligned}
\|\bU_n\|_{\hB^{0}_{3,r_2} (\BR^3)} 
\lesssim \big( n^{-\ep/6} \big)^3
\lesssim n^{-\ep/2},
\end{aligned}
\end{equation*}
which gives us that 
\begin{equation}\label{es:U_n}
\|\bU_n\|_{\hB^{0}_{3,r} (\BN(n))}  
 \lesssim \big(\sum_{k\in \BN(n)} 1^{r_1} \big)^{1/r_1}
 \|\bU_n\|_{\hB^{0}_{3,r_2} (\BN(n))}  
 \lesssim n^{-(\ep-2\sigma_2)/6}
\end{equation}
for $1/r_1=1/r-1/r_2=(\ep+\sigma_2)/3.$
Then Lemma \ref{lem:key} and \eqref{es:U_n} furnish that
\begin{equation*}
\begin{aligned}
\|\bu_n\|_{\hB^{0}_{3,r}(\BR^3)} \geq \|\bu_n\|_{\hB^{0}_{3,r}(\BN(n))} 
&\geq \|\bG_n\|_{\hB^{0}_{3,r}(\BN(n))} 
-\|\bg_n\|_{\hB^{0}_{3,r}(\BN(n))}
-\|\bU_n\|_{\hB^{0}_{3,r}(\BN(n))}\\
& \geq c-C_1 n^{-\ep/6}-C_2 n^{-(\ep-2\sigma_2)/6} \geq c/2
\end{aligned}
\end{equation*}
for $n$ large enough. This inflation phenomenon yields the discontinuity of the solution mapping from $\CD_{3}$ to $\CS_{3}.$

\section*{Acknowledgement}
J.T. is supported by the CY Initiative of Excellence, project CYNA (CY Nonlinear Analysis);
H.T. is supported by Grant-in-Aid for JSPS Fellows, Grant Number JP22KJ1642; 
X.Z. is partially supported by NSF of China under Grant 12101457 and  supported by the Fundamental Research Funds for the Central Universities.

Part of this work has been performed during a research visit of the second and third authors in Keio University. Prof. Linyu Peng is warmly thanked for his stimulating discussions during this period.


\end{document}